\documentclass[12pt,a4paper]{article}
\usepackage{amsfonts}
\usepackage{amsfonts,amssymb,amsmath,indentfirst,amsthm}

\input{tcilatex}

\begin{document}

\date{}
\title{Higher gradients estimates in Morrey
spaces for weak solutions to linear ultraparabolic equations}
\author{Yan Dong, Pengcheng Niu \thanks{%
Pengcheng Niu: Corresponding author, Department of Applied
Mathematics, Northwestern Polytechnical University, Xi'an, Shaanxi,
710129, China. e-mail: pengchengniu@nwpu.edu.cn}
}
\date{}
\maketitle

\newtheorem{thm}{{\indent}Theorem} \newtheorem{cor}[thm]{{\indent}Corollary}
\newtheorem{lem}[thm]{{\indent}Lemma} \newtheorem{prop}[thm]{%
{\indent}Proposition} \theoremstyle{definition} \newtheorem{defn}[thm]{%
{\indent}Definition}

\newtheorem{theo}{\hspace*{2.0em}Theorem\hspace*{0.1em}}[section]
\newenvironment{keywords}{\par\textbf{keywords:}\mbox{ }}{ } \newenvironment{%
coj}{\par\textbf{Conjecture:}\mbox{ }}{ } \numberwithin{equation}{section}

\textbf{Abstract.} The aim of this paper is to consider the linear
ultraparabolic equation with bounded and VMO coefficients $a_{ij}
(z)$. Assume that the operator $L_0$ obtained by freezing the
coefficients $a_{ij}(z)$ at any point ${z_0} \in {\mathbb{R}^{N +
1}}$ is hypoelliptic. We first establish a Caccioppoli type
inequality by choosing a cutoff function, a Sobolev type inequality
by prosperities of the fundamental solution to $L_0$, and a
Poincar\'{e} type inequality with a new cutoff function. Then $L^p$
estimate for weak solutions is derived by using the reverse
H\"{o}lder inequality on homogeneous spaces. Finally, higher Morrey
estimates for weak solutions to the above equation are shown by
investigating a homogeneous ultraparabolic equation of variable
coefficients with a nonhomogeneous boundary value condition, and a
nonhomogeneous ultraparabolic equation of variable coefficients with
homogeneous boundary value condition.




\section{Introduction}

\label{introduction}

In the paper, we consider the ultraparabolic equation of the kind
\begin{equation}
\label{eq1} Lu = \sum\limits_{i,j = 1}^{{m_0}} {{\partial _{{x_i}}}\left( {{a_{ij}}(z){\partial _{{x_j}}}u(z)} \right)}  + \sum\limits_{i,j = 1}^N {{b_{ij}}{x_i}{\partial _{{x_j}}}u(z)}  - {\partial _t}u(z) = g(z) + \sum\limits_{j = 1}^{{m_0}} {{\partial _{{x_j}}}{f_j}(z)} ,
\end{equation}
\noindent where $z = (x,t) \in {\mathbb{R}^{N + 1}}$, $1 \le {m_0}
\le N$, ${b_{ij}} \in \mathbb{R}$ ($i,j = 1, \ldots ,N$), $g,{f_j}
\in {L^p}\left( \Omega  \right)$ or ${L^{p,\lambda }}\left( \Omega  \right)$, ${L^{p,\lambda }}\left(
\Omega  \right)$ is a Morrey space, here $\Omega$ is a bounded
domain in $\mathbb{R}^{N + 1}$, $p \ge 2$, $0 \le \lambda  < Q +
2$, see Section 2 for the meaning of $Q$.

The assumptions to (\ref{eq1}) are

(H1) (ellipticity condition on $\mathbb{R}^{m_0}$) Let coefficients ${a_{ij}}(z) \in VMO
\cap {L^\infty }(\Omega )$(see Section 2 for the definition of VMO),
${a_{ij}}(z) = {a_{ji}}(z)$, satisfying that there exists a constant $\Lambda  > 1$ such that for any $z
\in {\mathbb{R}^{N + 1}}$, $\xi  \in {\mathbb{R}^{{m_0}}}$,
\[
{\Lambda ^{ - 1}}{\left| \xi  \right|^2} \le \sum\limits_{i,j =
1}^{{m_0}} {{a_{ij}}(z){\xi _i}{\xi _j}}  \le \Lambda {\left| \xi
\right|^2}.
\]

(H2) The constant matrix $B = {\left( {{b_{ij}}} \right)_{i,j = 1,
\ldots ,N}}$ in (\ref{eq1}) has the form
\[B = \left( {\begin{array}{*{20}{c}}
   0 & {{B_1}} & 0 &  \cdots  & 0  \\
   0 & 0 & {{B_2}} &  \cdots  & 0  \\
    \vdots  &  \vdots  &  \vdots  &  \ddots  &  \vdots   \\
   0 & 0 & 0 &  \cdots  & {{B_r}}  \\
   0 & 0 & 0 &  \cdots  & 0  \\
\end{array}} \right),\]
\noindent where ${B_k}(k = 1,2, \ldots ,r)$ is a ${m_{k - 1}}
\times {m_k}$ matrix with rank $m_k$ and
\[{m_0} \ge {m_1} \ge  \cdots  \ge {m_r} \ge 1, \quad \sum\limits_{k = 0}^r {{m_k}}  = N.\]

The equation (\ref{eq1}) can be written as
\[Lu = div\left( {A(z){D_0}u} \right) + Yu = g + divf,\]
\noindent where ${D_0} = \left( {{\partial _{{x_1}}},{\partial
_{{x_2}}}, \ldots ,{\partial _{{x_{{m_0}}}}},0, \ldots ,0} \right)$,
$Yu = \left\langle {x,BDu} \right\rangle  - {\partial _t}u$, $D =
\left( {{\partial _{{x_1}}},{\partial _{{x_2}}}, \ldots ,{\partial
_{{x_N}}}} \right)$, $f = \left( {{f_1},{f_2}, \ldots
,{f_{{m_0}}},0, \ldots ,0} \right)$, $A(z)$ is a $N \times
N$ matrix, namely,
\[A(z) = \left( {\begin{array}{*{20}{c}}
   {{A_0}(z)} & 0  \\
   0 & 0  \\
\end{array}} \right), \quad {A_0}(z) = {\left( {{a_{ij}}(z)} \right)_{i,j = 1, \ldots ,{m_0}}}.\]

Regularity for weak solutions to parabolic equations were provided
by many authors including DiBenedetto [5], Friendman [9], Krylov
[15], Ladyzhenskaya, Solonnikov, Ural'tseva [16], Lieberman [18] and
references therein.

In recent decades, many scholars have concerned with regularity of weak solutions to ultraparabolic equations.
These equations are closely related to finance, Brown motion, partical physics and human vision, etc.
The classic linear parabolic equation is usually of the form
\[\sum\limits_{i = 1}^N {{\partial _{{x_i}{x_j}}}u(x,t)}  - {\partial _t}u(x,t) = f(x,t).\]
But we see that (1.1) is strongly degenerate if $1 \le {m_0} < N$ and there is a drift $Yu$. These make research on regularity to (1.1) different from parabolic equation.

For the homogeneous ultraparabolic
equation
\begin{equation}
\label{eq2}  Lu = \sum\limits_{i,j = 1}^{{m_0}} {{\partial
_{{x_i}}}\left( {{a_{ij}}(z){\partial _{{x_j}}}u(z)} \right)}  +
\sum\limits_{i,j = 1}^N {{b_{ij}}{x_i}{\partial _{{x_j}}}u(z)}  -
{\partial _t}u(z) = 0,
\end{equation}
Polidoro in [23] got global lower bound of the fundamental solution
to (\ref{eq2}). The boundedness of weak solutions to (\ref{eq2})
with measurable coefficients was investigated by Pascucci and
Polidoro in [22] with Moser's iteration method based on a
combination of a Caccioppoli type estimate and the classical
embedding Sobolve inequality. Wang and Zhang in [25] obtained
H\"{o}lder estimates for weak solutions to (\ref{eq2}) with
measurable coefficients by deriving local a priori estimate to
(\ref{eq2}) and a Poincar\'{e} inequality of nonnegative weak lower
solution.

To the following ultraparabolic equation
\begin{equation}
\label{eq3} Lu = \sum\limits_{i,j = 1}^{{m_0}} {{\partial
_{{x_i}}}\left( {{a_{ij}}(z){\partial _{{x_j}}}u(z)} \right)}  +
\sum\limits_{i,j = 1}^N {{b_{ij}}{x_i}{\partial _{{x_j}}}u(z)}  -
{\partial _t}u(z) = \sum\limits_{j = 1}^{{m_0}} {{\partial
_{{x_j}}}{F_j}(x,t)},
\end{equation}
\noindent where ${F_j} \in L_{loc}^p\left( {{\mathbb{R}^{N + 1}}}
\right)(1 < p < \infty)$, coefficients $a_{ij}(z)$ belong to VMO
spaces, Manfredini and Polidoro in [19] established $L^p$ estimates
and H\"{o}lder continuity for weak solutions $u \in L_{loc}^p\left(
{{\mathbb{R}^{N + 1}}} \right)$. If ${F_j} \in L_{loc}^{p,\lambda
}\left( {{\mathbb{R}^{N + 1}}} \right)(1 < p < \infty, 0 \le \lambda
< Q + 2)$ and coefficients $a_{ij}(z)$ belong to some VMO spaces,
Polidoro and Ragusa in [24] derived H\"{o}lder regularity for weak
solution $u \in L_{loc}^p\left( {{\mathbb{R}^{N + 1}}} \right)$ to
(\ref{eq3}). Bramanti, Cerutti and Manfredini [1] proved local $L^p$
estimates for second order derivatives ${\partial
_{{x_i}{x_j}}}u\left( {i,j = 1, \ldots ,{m_0}} \right)$ of strong
solutions to the nondivergence ultraparabolic equation
\[\sum\limits_{i,j = 1}^{{m_0}} {{a_{ij}}(z){\partial _{{x_i}{x_j}}}u}  + \left\langle {x,BD} \right\rangle u - {\partial _t}u = f\]
with $a_{ij}(z)$ being in VMO and $f \in L^p$. The methods in [1,
19, 24] are based on the representation formulae for solutions and
estimates of singular integral operators. More related results also
see  Cinti, Passcucci and Polidoro [4], Xin and Zhang [26], Zhang
[27] and references therein.

The aim of this paper is to establish higher integrability for weak
solution $u \in W_2^{1,1}\left( \Omega  \right)$ to (\ref{eq1}) with
the method of a priori estimates. Results on higher integrability of
parabolic equations see Byun and Wang [2], Fugazzda [10], Palagachev
and Softova [21] and references therein. The first result here is
the higher $L^p$($p > 2$) estimate. For this purpose, an appropriate
frame is homogeneous spaces. Bramanti, Cerutti and Manfredini [1]
pointed out that the ball related to a quasidistance (see Section 2
below) is a homogeneous space and Gianazza [11] proved a reverse
H\"{o}lder inequality on homogeneous spaces. These facts will play
important roles. In spite of this, some new preliminary conclusions
are needed. Inspired by the way in [22], we deduce a Caccioppoli
type inequality and a Sobolev type inequality for weak solution to
(\ref{eq1}). Following to [7] and constructing suitable cutoff
functions, a Poincar\'{e} type inequality for weak solution to
(\ref{eq1}) is obtained. And then we prove higher $L^p$ estimates
for gradients of weak solutions to (\ref{eq1}) by using these new
inequalities and the reverse H\"{o}lder inequality on the
homogeneous spaces in [11].

The second result is on higher integrability in Morrey spaces for
gradients of weak solution $u \in W_2^{1,1}\left( \Omega \right)$ to
(\ref{eq1}). With the aid of the approach in the studying of
parabolic equations (e.g., see [12]), we consider a homogeneous
ultraparabolic equation of variable coefficients with a
nonhomogeneous boundary value condition, i.e., (6.1) below, and a
nonhomogeneous ultraparabolic equation of variable coefficients with
homogeneous boundary value condition, i.e., (6.2) below. The $L^p$
estimate for gradients of weak solutions to (6.1) is obtained by
proving a local ${L^\infty }$ estimate and a local $L^2$ estimate of
homogeneous ultraparabolic equation of constant coefficient, (5.1)
below. Then we establish a local $L^p$ estimate for gradients of
weak solutions to (6.2). These results are of independent interest.
Finally, higher integrability in Morrey spaces for gradients of weak
solutions to (\ref{eq1}) is deduced by using a known iteration
lemma.

The following is the notion of weak solution to (1.1).

\textbf{Definition 1.1} If $u \in W_2^{1,1}\left( \Omega \right)$
and for any $\psi  \in C_0^\infty \left( \Omega  \right)$,
\[ - \int_\Omega  {A{D_0}u{D_0}\psi } dz + \int_\Omega  {\psi Yu} dz = \int_\Omega  {\left( {g\psi  - f{D_0}\psi } \right)} dz,\]
\noindent then we say that $u$ is a weak solution to (\ref{eq1}).

The main results of this paper are stated as follows.

\textbf{Theorem 1.1} Suppose that assumptions (H1) and (H2) hold. If $u
\in W_2^{1,1}\left( \Omega  \right)$ is a weak solution to
(\ref{eq1}), $g,{f_j} \in {L^p}\left( \Omega  \right)$, then there exists a constant ${\varepsilon _0} > 0$
such that for any $p \in \left[ {2,2 + \frac{{2Q}}{{Q +
2}}{\varepsilon _0}} \right)$, we have ${D_0}u \in L_{loc}^p(\Omega
)$ and for any $\Omega ' \subset  \subset \Omega '' \subset \subset
\Omega$,
\begin{equation}
\label{eq4} {\left\| {{D_0}u} \right\|_{{L^p}(\Omega ')}} \le
c\left( {{{\left\| {{D_0}u} \right\|}_{{L^2}(\Omega '')}} +
{{\left\| g \right\|}_{{L^p}(\Omega )}} + {{\left\| f
\right\|}_{{L^p}(\Omega )}}} \right).
\end{equation}

\textbf{Theorem 1.2} Under (H1) and (H2), let $u \in W_2^{1,1}\left(
\Omega  \right)$ be a weak solution to (\ref{eq1}), $g,{f_j} \in
{L^{p,\lambda }}\left( \Omega  \right)$, then for any $p \in \left[
{2,2 + \frac{{2Q}}{{Q + 2}}{\varepsilon _0}} \right)$, ${\varepsilon
_0}$ as in Theorem 1.1, we have ${D_0}u \in L_{loc}^{p,\lambda
}(\Omega )$ ($0 < \lambda  < Q + 2)$ and for any $\Omega ' \subset
\subset \Omega '' \subset \subset \Omega$,
\begin{equation}
\label{eq5} {\left\| {{D_0}u} \right\|_{{L^{p,\lambda }}(\Omega ')}}
\le c\left( {{{\left\| {{D_0}u} \right\|}_{{L^2}(\Omega '')}} +
{{\left\| g \right\|}_{{L^{p,\lambda }}(\Omega )}} + {{\left\| f
\right\|}_{{L^{p,\lambda }}(\Omega )}}} \right).
\end{equation}
\noindent

This paper is organized as follows. In Section 2, we describe some
basic knowledge and some known material on the frozen operator $L_0$
of $L$ and the fundamental solution of $L_0$, and collect several
useful lemmas which will be used later on. Section 3 is devoted to
proofs of a Caccioppoli type inequality, a Sobolev type inequality
and a Poincar\'{e} type inequality for weak solutions. In Section 4,
the proof of Theorem 1.1 is given by using the inequalities in
Section 3 and the reverse H\"{o}lder inequality in [11]. In Section
5, we derive a higher $L^p$ estimate for gradient of weak solutions
to (5.1). In Section 6, the proof of Theorem 1.2 is ended by local
$L^p$ estimate for gradient of weak solutions to (6.1) and (6.2).


\section{Preliminaries}

\label{2}

For any ${z_0} \in \Omega  \subset {\mathbb{R}^{N + 1}}$, we denote
the frozen operator of $L$ by
\begin{equation}
\label{eq6} {L_0} = \sum\limits_{i,j = 1}^{{m_0}} {{\partial
_{{x_i}}}\left( {{a_{ij}}({z_0}){\partial _{{x_j}}}} \right)}  +
\sum\limits_{i,j = 1}^N {{b_{ij}}{x_i}{\partial _{{x_j}}}}  -
{\partial _t}.
\end{equation}
Now one can introduce the following.

\textbf{Definition 2.1} For any $(x,t),(\xi ,\tau ) \in
{\mathbb{R}^{N + 1}}$, set a multiplication law
\[
(x,t) \circ (\xi ,\tau ) = \left( {\xi  + E(\tau )x,t + \tau }
\right), \quad E(\tau ) = \exp ( - \tau {B^T}).
\]
We say that $\left( {{\mathbb{R}^{N + 1}}, \circ } \right)$ is a
noncommutative Lie group with neutral element $(0,0)$, the inverse
of an element $(x,t) \in {\mathbb{R}^{N + 1}}$ is
\[{(x,t)^{ - 1}} = \left( { - E( - t)x, - t} \right).\]

Authors in [17] claimed that the frozen operator $L_0$ is
hypoelliptic and left invariant about the groups of translations and
dilations. In this case, the dilations associated to $L_0$ are given
by
\[{\delta _\lambda } = diag\left( {\lambda {I_{{m_0}}},{\lambda ^3}{I_{{m_1}}},
 \ldots ,{\lambda ^{2r + 1}}{I_{{m_r}}},{\lambda ^2}} \right), \quad \lambda  > 0,\]
\noindent here ${I_{{m_k}}}$ denotes the ${m_k} \times {m_k}$
identity matrix, and
\[\det \left( {{\delta _\lambda }} \right) = {\lambda ^{Q + 2}},\]
\noindent with $Q + 2 = {m_0} + 3{m_1} +  \cdots  + (2r + 1){m_r} +
2$. The number $Q+2$ is called the homogeneous dimension of $\mathbb{R}^{N
+ 1}$, and $Q$ the homogeneous dimension of $\mathbb{R}^N$. Note
that $L_0$ is ${\delta _\lambda }$ homogeneous of degree 2, namely,
for any $\lambda  > 0$,
\[{L_0} \circ {\delta _\lambda } = {\lambda ^2}\left( {{\delta _\lambda } \circ {L_0}} \right).\]

Due to [14], the fundamental solution ${\Gamma _0}( \cdot ,\zeta )$
of $L_0$ has an explicit expression in the pole $\zeta  \in
{\mathbb{R}^{N + 1}}$, which is, for any $z, \zeta  \in
{\mathbb{R}^{N + 1}}$, $z \ne \zeta$,
\begin{equation}
\label{eq7} {\Gamma _0}(z,\zeta ) = {\Gamma _0}({\zeta ^{ - 1}}
\circ z,0),
\end{equation}
\noindent where
\[{\Gamma _0}\left( {(x,t),(0,0)} \right) = \left\{ \begin{array}{l}
 \frac{1}{{{{\left( {{{\left( {4\pi } \right)}^N}\det C(t)} \right)}^{\frac{1}{2}}}}}\exp \left( { - \frac{1}{4}\left\langle {{C^{ - 1}}(t)x,x} \right\rangle } \right),t > 0, \\
 0,{\kern 1pt} {\kern 1pt} {\kern 1pt} {\kern 1pt} {\kern 1pt} {\kern 1pt} {\kern 1pt} {\kern 1pt} {\kern 1pt} {\kern 1pt} {\kern 1pt} {\kern 1pt} {\kern 1pt} {\kern 1pt} {\kern 1pt} {\kern 1pt} {\kern 1pt} {\kern 1pt} {\kern 1pt} {\kern 1pt} {\kern 1pt} {\kern 1pt} {\kern 1pt} {\kern 1pt} {\kern 1pt} {\kern 1pt} {\kern 1pt} {\kern 1pt} {\kern 1pt} {\kern 1pt} {\kern 1pt} {\kern 1pt} {\kern 1pt} {\kern 1pt} {\kern 1pt} {\kern 1pt} {\kern 1pt} {\kern 1pt} {\kern 1pt} {\kern 1pt} {\kern 1pt} {\kern 1pt} {\kern 1pt} {\kern 1pt} {\kern 1pt} {\kern 1pt} {\kern 1pt} {\kern 1pt} {\kern 1pt} {\kern 1pt} {\kern 1pt} {\kern 1pt} {\kern 1pt} {\kern 1pt} {\kern 1pt} {\kern 1pt} {\kern 1pt} {\kern 1pt} {\kern 1pt} {\kern 1pt} {\kern 1pt} {\kern 1pt} {\kern 1pt} {\kern 1pt} {\kern 1pt} {\kern 1pt} {\kern 1pt} {\kern 1pt} {\kern 1pt} {\kern 1pt} {\kern 1pt} {\kern 1pt} {\kern 1pt} {\kern 1pt} {\kern 1pt} {\kern 1pt} {\kern 1pt} {\kern 1pt} {\kern 1pt} {\kern 1pt} {\kern 1pt} {\kern 1pt} {\kern 1pt} {\kern 1pt} {\kern 1pt} {\kern 1pt} {\kern 1pt} {\kern 1pt} {\kern 1pt} {\kern 1pt} {\kern 1pt} {\kern 1pt} {\kern 1pt} {\kern 1pt} {\kern 1pt} {\kern 1pt} {\kern 1pt} {\kern 1pt} {\kern 1pt} {\kern 1pt} {\kern 1pt} {\kern 1pt} {\kern 1pt} {\kern 1pt} {\kern 1pt} {\kern 1pt} {\kern 1pt} {\kern 1pt} {\kern 1pt} {\kern 1pt} {\kern 1pt} {\kern 1pt} {\kern 1pt} {\kern 1pt} {\kern 1pt} {\kern 1pt} {\kern 1pt} {\kern 1pt} {\kern 1pt} {\kern 1pt} {\kern 1pt} {\kern 1pt} {\kern 1pt} {\kern 1pt} {\kern 1pt} {\kern 1pt} {\kern 1pt} {\kern 1pt} {\kern 1pt} {\kern 1pt} {\kern 1pt} {\kern 1pt} {\kern 1pt} {\kern 1pt} {\kern 1pt} {\kern 1pt} {\kern 1pt} {\kern 1pt} {\kern 1pt} {\kern 1pt} {\kern 1pt} {\kern 1pt} {\kern 1pt} {\kern 1pt} {\kern 1pt} {\kern 1pt} {\kern 1pt} {\kern 1pt} {\kern 1pt} {\kern 1pt} {\kern 1pt} {\kern 1pt} {\kern 1pt} {\kern 1pt} {\kern 1pt} {\kern 1pt} {\kern 1pt} {\kern 1pt} {\kern 1pt} {\kern 1pt} {\kern 1pt} {\kern 1pt} {\kern 1pt} {\kern 1pt} {\kern 1pt} {\kern 1pt} {\kern 1pt} {\kern 1pt} {\kern 1pt} {\kern 1pt} {\kern 1pt} {\kern 1pt} {\kern 1pt} {\kern 1pt} {\kern 1pt} {\kern 1pt} {\kern 1pt} {\kern 1pt} {\kern 1pt} {\kern 1pt} {\kern 1pt} {\kern 1pt} {\kern 1pt} {\kern 1pt} {\kern 1pt} {\kern 1pt} t \le 0, \\
 \end{array} \right.\]
\[C(t) = \int_0^t {E(s){A_0}{E^T}(s)} ds.\]

It is known that $C(t)$ is strictly
positive for every positive $t$. In view of the invariance
properties of $L_0$, we have that for any $z \in {\mathbb{R}^{N +
1}}\backslash \{ 0\}$ and $\lambda  > 0$,
\[{\Gamma _0}\left( {{\delta _\lambda }(z),0} \right) = {\lambda ^{ - Q}}{\Gamma _0}(z,0).\]
We also observe that ${\Gamma _0}$ is ${\delta _\lambda }$
homogeneous of degree $-Q$. For $i,j = 1, \ldots ,{m_0}$,
${D_{{x_i}}}{\Gamma _0}$ and ${D_{{x_i}{x_j}}}{\Gamma _0}$ are
${\delta _\lambda }$ homogeneous of degree $-(Q+1)$ and $-(Q+2)$,
respectively.

For any $(x,t) \in {\mathbb{R}^{N + 1}}$, the homogeneous
norm of $(x,t)$ with respect to ${\delta _\lambda }$ is defined by
\[\left\| {(x,t)} \right\| = \sum\limits_{j = 1}^N {{{\left| {{x_j}} \right|}^{\frac{1}{{{\alpha _j}}}}}}  + {\left| t \right|^{\frac{1}{2}}},\]
\noindent where ${\alpha _j} = 1$, if $1 \le j \le {m_0}$; ${\alpha
_j} = 2k + 3$, if ${m_k} < j \le {m_{k + 1}}\left( {0 \le k \le r -
1} \right)$. For any $z,\zeta  \in {\mathbb{R}^{N + 1}}$, we denote
the quasidistance by
\[d\left( {z,\zeta } \right) = \left\| {{\zeta ^{ - 1}} \circ z} \right\|.\]

\textbf{Lemma 2.2 }([6, Lemma 2.1]) For any bounded domain $\Omega
\subset {\mathbb{R}^{N + 1}}$, $d\left( {z,\zeta } \right)$ is a
quasisymmetric quasidistance in $\Omega$, if for any $z,z',\zeta \in
{\mathbb{R}^{N + 1}}$,
\[d\left( {z,\zeta } \right) \le cd\left( {\zeta ,z} \right),\quad
d\left( {z,\zeta } \right) \le c\left( {d\left( {\zeta ,z'} \right) + d\left( {z',\zeta } \right)} \right).\]

The ball with respect to $d$ centered at $z_0$ is denoted by
\[{B_R}({z_0}) = B({z_0},R) = \left\{ {\zeta  \in {\mathbb{R}^{N + 1}}:{\kern 1pt} {\kern 1pt} {\kern 1pt} d\left( {{z_0},\zeta } \right) < R} \right\}.\]
Note clearly that $B(0,R) = {\delta _R}B(0,1)$.

\textbf{Remark 2.3 }Recalling [1, Remark 1.5], it holds that for any
${z_0} \in {\mathbb{R}^{N + 1}}$, $R > 0$,
\[\left| {B({z_0},R)} \right| = \left| {B(0,R)} \right| = \left| {B(0,1)} \right|{R^{Q + 2}},\]
\[\left| {B({z_0},2R)} \right| = {2^{Q + 2}}\left| {B({z_0},R)} \right|,\]
\noindent and therefore the space $\left( {{\mathbb{R}^{N +
1}},dz,d} \right)$ is a homogeneous space. The fact allows us to
employ known conclusions in homogeneous spaces.

If one does not need to concern the center of the ball,
$B({z_0},R)$ can simply be written as $B_R$. For convenience, we usually
consider the estimates on cubes instead of balls. Let us describe the notion of cubes. For any $x =
(x',\bar x)$, $x' = \left( {{x_1}, \cdots ,{x_{{m_0}}}} \right)$,
$\bar x = \left( {{x_{{m_0} + 1}}, \cdots ,{x_N}} \right)$, the cube is denoted by
\[\begin{array}{l}
 {Q_R} = \left\{ {(x,t)\left| {{t_0} - {R^2}/2 \le t \le {t_0} + {R^2}/2,\left| {x'} \right| \le R,\left| {{x_{{m_0} + 1}}} \right|
 \le {{\left( {\Lambda {N^2}R} \right)}^3}, \cdots ,} \right.} \right. \\
 \quad \quad \quad \left. {\left| {{x_N}} \right| \le {{\left( {\Lambda {N^2}R} \right)}^{2r + 1}}} \right\} \\
 \end{array}\]
Also, we write
\[{I_R} = \left[ {{t_0} - \frac{{{R^2}}}{2},{t_0} + \frac{{{R^2}}}{2}} \right],\]
\[{K_R} = \left\{ {x'\left| {\left| {x'} \right| \le R} \right.} \right\},\]
\[{S_R} = \left\{ {\bar x\left| {\left| {{x_{{m_0} + 1}}} \right|
\le {{\left( {\Lambda {N^2}R} \right)}^3}, \cdots ,\left| {{x_N}}
\right| \le {{\left( {\Lambda {N^2}R} \right)}^{2r + 1}}} \right.}
\right\},\] \noindent then ${Q_R} = {K_R} \times {S_R} \times
{I_R}$.

A cube of centered at $(0,0)$ is simply denoted by
\[{Q_R}(0,0) = \left\{ {(x,t)\left| {\left| t \right| \le {R^2},\left| {{x_1}} \right| \le {R^{{\alpha _1}}}, \cdots ,\left| {{x_N}} \right| \le {R^{{\alpha _N}}}} \right.} \right\}.\]
It is easy to find that there exists a constant ${c_0} = {c_0}\left(
{B,N} \right) > 0$, such that
\[{Q_{R/{c_0}}}(0,0) \subset {B_R}(0,0) \subset {Q_{{c_0}R}}(0,0).\]
We state a result on ${\delta _\lambda }$ homogeneous functions in
[8, 22].

\textbf{Lemma 2.4} Let $\alpha  \in \left[ {0,Q + 2} \right]$ and $G
\in C\left( {{\mathbb{R}^{N + 1}}\backslash \{ 0\} } \right)$ be a
${\delta _\lambda }$ homogeneous function of degree $\alpha  - Q -
2$. If $f \in {L^p}\left( {{\mathbb{R}^{N + 1}}} \right)$, $p \in
\left[ {1, + \infty } \right)$, then the function
$${G_f}(z) \equiv \int_{{\mathbb{R}^{N + 1}}} {G({\zeta ^{
- 1}} \circ z)f(\zeta )} d\zeta $$
\noindent is defined almost
everywhere and there exists a constant $c = c\left( {Q,P} \right) > 0$
such that
\begin{equation}
\label{eq8} {\left\| {{G_f}} \right\|_{{L^q}({\mathbb{R}^{N + 1}})}}
\le c\mathop {\max }\limits_{\left\| z \right\| = 1} \left| {G(z)}
\right|{\left\| f \right\|_{{L^p}({\mathbb{R}^{N + 1}})}},
\end{equation}
\noindent where $\frac{1}{q} = \frac{1}{p} - \frac{\alpha }{{Q +
2}}$.

This lemma can be used to yield the following.

\textbf{Lemma 2.5 }Let $f \in {L^{\frac{{2(Q + 2)}}{{Q + 4}}}}\left(
{{\mathbb{R}^{N + 1}}} \right)$. There exists a positive constant $c
= c(Q)$ such that
\begin{equation}
\label{eq9} {\left\| {{\Gamma _0}(f)} \right\|_{{L^{\frac{{2(Q +
2)}}{Q}}}({\mathbb{R}^{N + 1}})}} \le c{\left\| f
\right\|_{{L^{\frac{{2(Q + 2)}}{{Q + 4}}}}({\mathbb{R}^{N + 1}})}},
\end{equation}
\begin{equation}
\label{eq10} {\left\| {{\Gamma _0}({D_0}f)}
\right\|_{{L^2}({\mathbb{R}^{N + 1}})}} \le c{\left\| f
\right\|_{{L^{\frac{{2(Q + 2)}}{{Q + 4}}}}({\mathbb{R}^{N + 1}})}},
\end{equation}
where ${\Gamma _0}\left( f \right)(z) = \int_{{\mathbb{R}^{N + 1}}} {{\Gamma _0}\left( {z,\zeta } \right)f\left( \zeta  \right)} d\zeta $, ${\Gamma _0}\left( {{D_0}f} \right)(z) = \int_{{\mathbb{R}^{N + 1}}} {{\Gamma _0}\left( {z,\zeta } \right){D_0}f\left( \zeta  \right)} d\zeta$.

\textbf{Proof:} Since ${\Gamma _0}$ is homogeneous of degree $-Q$
with respect to ${\delta _\lambda }$, we immediately have
(\ref{eq9}) from Lemma 2.4 by taking $\alpha  = 2$, $q = \frac{{2(Q
+ 2)}}{Q}$ and $p = \frac{{2(Q + 2)}}{{Q + 4}}$. Noting that ${\partial _{{x_i}}}{\Gamma _0}$ is homogeneous of
degree $ - \left( {Q + 1} \right)$ with respect to ${\delta _\lambda
}$, (\ref{eq10}) holds by Lemma 2.4 with $\alpha  = 1$, $q=2$ and
$p = \frac{{2(Q + 2)}}{{Q + 4}}$.

\textbf{Definition 2.6 }(Morrey space ${L^{p,\lambda }}$) Let
$\Omega$ be an open subset in ${\mathbb{R}^{N + 1}}$, $1 \le p <  +
\infty $, $\lambda
> 0$. We say that $f \in {L^p}(\Omega )$ belongs to the Morrey space ${L^{p,\lambda }}(\Omega )$, if
\[{\left\| f \right\|_{{L^{p,\lambda }}}} = \mathop {\sup }\limits_{{z_0} \in \Omega ,\rho  > 0}
{\left( {\frac{{{\rho ^\lambda }}}{{\left| {\Omega  \cap {B_\rho
}({z_0})} \right|}} \int_{\Omega  \cap {B_\rho }({z_0})} {{{\left| f
\right|}^p}} dz} \right)^{\frac{1}{p}}} < \infty .\]

\textbf{Definition 2.7} (Sobolev space $W_p^{1,1}$) Let $\Omega$ be
an open subset in ${\mathbb{R}^{N + 1}}$. The Sobolev space with respect to (1.1) is
defined by $$W_p^{1,1}\left( \Omega  \right) = \left\{ {u \in
{L^p}(\Omega ):{\partial _{{x_i}}}u,Yu \in {L^p}(\Omega ),{\kern
1pt} {\kern 1pt} {\kern 1pt} i,j = 1, \ldots ,{m_0}} \right\}$$
\noindent with the norm $$\left\| u \right\|_{W_p^{1,1}\left( \Omega
\right)}^p = \left\| u \right\|_{{L^p}\left( \Omega  \right)}^p +
\sum\limits_{i = 1}^{{m_0}} {\left\| {{\partial _{{x_i}}}u}
\right\|_{{L^p}\left( \Omega  \right)}^p}  + \left\| {Yu}
\right\|_{{L^p}\left( \Omega  \right)}^p.$$

The space $W_{p,0}^{1,1}\left( \Omega  \right)$ is the closure of
$C_0^\infty \left( {\bar \Omega } \right)$ in $W_p^{1,1}\left(
\Omega  \right)$.

\textbf{Definition 2.8 }(BMO and VMO spaces) For any $a \in
L_{loc}^1(\Omega )$, we set
\[{\eta _R}\left( a \right) =  \mathop {\sup }\limits_{{z_0} \in \Omega ,
0 \le \rho  \le R} \left( {\frac{1}{{\left| {\Omega  \cap {B_\rho
}({z_0})} \right|}}\int_{\Omega  \cap {B_\rho }({z_0})} {\left|
{a(z) - {a_{\Omega  \cap {B_\rho }({z_0})}}(z)} \right|} dz}
\right),\] \noindent where ${a_{\Omega  \cap {B_\rho }({z_0})}} =
\frac{1}{{\left| {\Omega  \cap {B_\rho }({z_0})}
\right|}}\int_{\Omega  \cap {B_\rho }({z_0})} {a(z)} dz$. If
$\mathop {\sup }\limits_{R > 0} {\eta _R}\left( a \right) < \infty
$, we say $a \in BMO(\Omega )$ (Bounded Mean Oscillation). Moreover,
if ${\eta _R}\left( a \right) \to 0$ as $R \to 0$, we say $a \in
VMO(\Omega )$ (Vanishing Mean Oscillation).

It is stated two iteration lemmas.

\textbf{Lemma 2.9} ([3]) Let $\varphi (t)$ be a bounded nonnegative
function on $[{T_0},{T_1}]$, where ${T_1} > {T_0} \ge 0$. Suppose
that for any $s,t:{T_0} \le t < s \le {T_1}$, $\varphi$ satisfies
\[\varphi (t) \le {\theta _1}\varphi (s) + \frac{{{a_2}}}{{{{(s - t)}^\alpha }}} + {b_2},\]
\noindent where ${\theta _1},{a_2},{b_2}$ and $\alpha $ are
nonnegative constants, and ${\theta _1} < 1$. Then for any ${T_0}
\le \rho  < R \le {T_1}$,
\begin{equation}
\label{eq11} \varphi (\rho ) \le c\left( {\frac{{{a_2}}}{{{{(R -
\rho )}^\alpha }}} + {b_2}} \right),
\end{equation}
\noindent where $c$ depends only on $\alpha$ and ${\theta _1}$.

\textbf{Lemma 2.10} (see [13, 20]) Let $H$ be a nonnegative
increasing function. Suppose that for any $\rho  < R \le {R_0} =
dist({z_0},\partial \Omega )$,
\[H(\rho ) \le {A_1}\left[ {{{\left( {\frac{\rho }{R}} \right)}^{{a_1}}} + \varepsilon } \right]H(R) + {B_1}{R^{{b_1}}},\]
\noindent where ${A_1},{a_1}$ and ${b_1}$ are positive constants
with ${a_1} > {b_1}$. Then there exist positive constants
${\varepsilon _1} = {\varepsilon _1}\left( {{A_1},{a_1},{b_1}}
\right)$ and $c = c\left( {{A_1},{a_1},{b_1}} \right)$, such that if
$\varepsilon  < {\varepsilon _1}$, then
\begin{equation}
\label{eq12} H(\rho ) \le c\left[ {{{\left( {\frac{\rho }{R}}
\right)}^{{b_1}}}H(R) + {B_1}{\rho ^{{b_1}}}} \right].
\end{equation}


\section{Preliminary inequalities}

\label{3}

\textbf{Theorem 3.1 }(Caccioppoli type inequality) Let $u \in
W_2^{1,1}\left( \Omega  \right)$ be a weak solution to (\ref{eq1}).
Then for any ${B_R} \subset \Omega$, $\rho  < R$, we have
\begin{equation}
\label{eq13} \int_{{B_\rho }} {{{\left| {{D_0}u} \right|}^2}} dz \le
\frac{c}{{{{(R - \rho )}^2}}}\int_{{B_R}} {{{\left| u \right|}^2}}
dz + c\int_{{B_R}} {\left( {{{\left| g \right|}^2} + {{\left| f
\right|}^2}} \right)} dz.
\end{equation}

\textbf{Proof:} Let $\xi (z) \in C_0^\infty ({B_R})$ be a cutoff
function satisfying
\begin{equation}
\label{eq14} \xi (z) = 1 (\left| z \right| < \rho), \xi (z) = 0
(\left| z \right| \ge R), 0 \le \xi  \le 1, \left| {{\partial
_{{x_j}}}\xi } \right|,\left| {{\partial _t}\xi } \right| \le
\frac{c}{{R - \rho }} (j = 1, \ldots ,N).
\end{equation}
Hence $$\left| {Y\xi } \right|{\rm{ = }}\left| {xBD\xi  - {\partial
_t}\xi } \right| \le c\left| {D\xi } \right|{\rm{ + }}c\left|
{{\partial _t}\xi } \right| \le \frac{c}{{R - \rho }},$$
and by the divergence theorem,
\[\int_{{B_R}} {Y\left( {{u^2}{\xi ^2}} \right)} dz = 0.\]
Multiplying both sides of (\ref{eq1}) by $u{\xi ^2}$ and integrating
on $B_R$, we have
\[\int_{{B_R}} {\left[ { - A{D_0}u{D_0}\left( {u{\xi ^2}} \right) + u{\xi ^2}Yu} \right]} dz
= \int_{{B_R}} {\left[ {gu{\xi ^2} - f{D_0}\left( {u{\xi ^2}}
\right)} \right]} dz\] and
\begin{align}
& \quad \int_{{B_R}} {A{\xi ^2}{D_0}u{D_0}u} dz \nonumber\\
& =  - 2\int_{{B_R}} {Au\xi {D_0}u{D_0}\xi } dz - \int_{{B_R}}
{{u^2}\xi Y\xi } dz - \int_{{B_R}} {gu{\xi ^2}} dz + \int_{{B_R}}
{f{\xi ^2}{D_0}u} dz \nonumber\\
& \quad + 2\int_{{B_R}} {fu\xi {D_0}\xi } dz.\label{3.3}%
\end{align}
By using (H1) and Young's inequality, it follows
\begin{align}
 & \quad {\Lambda ^{ - 1}}\int_{{B_R}} {{{\left| {{D_0}u} \right|}^2}{\xi ^2}} dz \nonumber\\
  & \le {c_\varepsilon }\int_{{B_R}} {{{\left| u \right|}^2}{{\left| {{D_0}\xi } \right|}^2}} dz
  + \varepsilon \int_{{B_R}} {{{\left| {{D_0}u} \right|}^2}{\xi ^2}} dz
  + \int_{{B_R}} {{{\left| u \right|}^2}\left| {Y\xi } \right|\xi } dz \nonumber\\
 &\quad + {c_\varepsilon }\int_{{B_R}} {{{\left| g \right|}^2}{\xi ^2}} dz
 + \varepsilon \int_{{B_R}} {{{\left| u \right|}^2}{\xi ^2}} dz
 + {c_\varepsilon }\int_{{B_R}} {{{\left| f \right|}^2}{\xi ^2}} dz
 + \varepsilon \int_{{B_R}} {{{\left| {{D_0}u} \right|}^2}{\xi ^2}} dz \nonumber\\
 &\quad + {c_\varepsilon }\int_{{B_R}} {{{\left| f \right|}^2}{\xi ^2}} dz
 + \varepsilon \int_{{B_R}} {{{\left| u \right|}^2}{{\left| {{D_0}\xi } \right|}^2}} dz \nonumber\\
 & \le \int_{{B_R}} {{{\left| u \right|}^2}\left( {{c_\varepsilon }{{\left| {{D_0}\xi } \right|}^2}
 + \left| {Y\xi } \right|\xi  + \varepsilon {\xi ^2} + \varepsilon {{\left| {{D_0}\xi } \right|}^2}} \right)} dz \nonumber\\
  &\quad + 2\varepsilon \int_{{B_R}} {{{\left| {{D_0}u} \right|}^2}{\xi ^2}} dz +
{c_\varepsilon }\int_{{B_R}} {\left( {{{\left| g \right|}^2} +
{{\left| f \right|}^2}} \right){\xi ^2}} dz.
\label{3.4}%
\end{align}
Choosing $\varepsilon$ small enough such that ${\Lambda ^{ - 1}} -
2\varepsilon  > 0$ and using the property of $\xi $, one has
\[\begin{array}{l}
 \int_{{B_R}} {{{\left| {{D_0}u} \right|}^2}{\xi ^2}} dz \\
  \le \int_{{B_R}} {{{\left| u \right|}^2}\left( {{c_\varepsilon }{{\left| {{D_0}\xi } \right|}^2} + \left| {Y\xi } \right|\xi  + \varepsilon {\xi ^2} + \varepsilon {{\left| {{D_0}\xi } \right|}^2}} \right)} dz + {c_\varepsilon }\int_{{B_R}} {\left( {{{\left| g \right|}^2} + {{\left| f \right|}^2}} \right){\xi ^2}} dz \\
  \le \int_{{B_R}} {{{\left| u \right|}^2}\left( {\frac{{{c_\varepsilon }}}{{{{(R - \rho )}^2}}} + \frac{{c\xi }}{{R - \rho }} + \varepsilon {\xi ^2} + \frac{\varepsilon }{{{{(R - \rho )}^2}}}} \right)} dz + {c_\varepsilon }\int_{{B_R}} {\left( {{{\left| g \right|}^2} + {{\left| f \right|}^2}} \right){\xi ^2}} dz. \\
 \end{array}\]
Consequently (\ref{eq13}) is proved.

\textbf{Theorem 3.2 }(Sobolev type inequality) Let $u \in W_2^{1,1}\left(
\Omega  \right)$ be a weak solution to (\ref{eq1}). Then for any
${B_R} \subset \Omega $, $\rho  < R$, it follows
\begin{align}
& \quad {\left\| u \right\|_{{L^2}\left( {{B_\rho }}
\right)}}\nonumber\\ & \le \frac{c}{{R - \rho }}\left( {{{\left\| u
\right\|}_{{L^{\frac{{2(Q + 2)}}{{Q + 4}}}}\left( {{B_R}} \right)}}
+ {{\left\| {{D_0}u} \right\|}_{{L^{\frac{{2(Q + 2)}}{{Q +
4}}}}\left( {{B_R}} \right)}} + {{\left\| g
\right\|}_{{L^{\frac{{2(Q + 2)}}{{Q + 4}}}}\left( {{B_R}} \right)}}
+ {{\left\| f \right\|}_{{L^{\frac{{2(Q + 2)}}{{Q + 4}}}}\left(
{{B_R}} \right)}}} \right).
\label{3.5}%
\end{align}

\textbf{Proof: }We represent $u$ in terms of the fundamental
solution ${\Gamma _0}$ of $L_0$ and apply the cutoff function $\xi$
in (3.2). For any $z \in {B_R}$,
\begin{equation}
\label{eq17} \left( {\xi u} \right)(z) = \int_{{B_R}} {\left[
{\left\langle {{A_0}{D_0}\left( {\xi u} \right),{D_0}{\Gamma _0}}
\right\rangle  - {\Gamma _0}Y\left( {\xi u} \right)} \right]} d\zeta
\buildrel \Delta \over = {I_1}(z) + {I_2}(z) + {I_3}(z),
\end{equation}
\noindent where $${I_1}(z) = \int_{{B_R}} {\left[ {{A_0}u{D_0}\xi
{D_0}{\Gamma _0} - {\Gamma _0}uY\xi } \right]} d\zeta ,$$

$${I_2}(z) = \int_{{B_R}} {\left[ {\left( {{A_0} - A} \right)\xi
{D_0}u{D_0}{\Gamma _0} - {\Gamma _0}A{D_0}u{D_0}\xi } \right]}
d\zeta $$
and
$${I_3}(z) = \int_{{B_R}} {\left[ {A{D_0}u{D_0}\left(
{\xi {\Gamma _0}} \right) - \xi {\Gamma _0}Yu} \right]} d\zeta .$$

It yields by using (\ref{eq9}) and (\ref{eq10}) that
\begin{align}
 &\quad {\left\| {{I_1}} \right\|_{{L^2}\left( {{B_R}} \right)}}
 \le 2{\left\| {\int_{{B_R}} {{A_0}u{D_0}\xi {D_0}{\Gamma _0}} d\zeta } \right\|_{{L^2}\left( {{B_R}} \right)}}
 + 2{\left\| {\int_{{B_R}} {{\Gamma _0}uY\xi } d\zeta } \right\|_{{L^2}\left( {{B_R}} \right)}} \nonumber\\
  & \le c{\left\| {{\Gamma _0}\left( {{D_0}\left( {u{D_0}\xi } \right)} \right)} \right\|_{{L^2}\left( {{B_R}} \right)}}
  + c{\left\| {{\Gamma _0}\left( {uY\xi } \right)} \right\|_{{L^2}\left( {{B_R}} \right)}} \nonumber\\
  & \le c{\left\| {u{D_0}\xi } \right\|_{{L^{\frac{{2(Q + 2)}}{{Q + 4}}}}\left( {{B_R}} \right)}}
  + c{\left\| {{\Gamma _0}\left( {uY\xi } \right)} \right\|_{{L^{\frac{{2(Q + 2)}}{Q}}}\left( {{B_R}} \right)}}{\left| {{B_R}} \right|^{\frac{1}{{Q + 2}}}} \nonumber\\
  &\le c{\left\| {u{D_0}\xi } \right\|_{{L^{\frac{{2(Q + 2)}}{{Q + 4}}}}\left( {{B_R}} \right)}}
  + c{\left\| {uY\xi } \right\|_{{L^{\frac{{2(Q + 2)}}{{Q + 4}}}}\left( {{B_R}} \right)}}R \nonumber\\
& \le \frac{c}{{R - \rho }}{\left\| u \right\|_{{L^{\frac{{2(Q +
2)}}{{Q + 4}}}}\left( {{B_R}} \right)}} + \frac{{cR}}{{R - \rho
}}{\left\| u \right\|_{{L^{\frac{{2(Q + 2)}}{{Q + 4}}}}\left(
{{B_R}} \right)}} \nonumber\\
& \le \frac{c}{{R - \rho }}{\left\| u
\right\|_{{L^{\frac{{2(Q + 2)}}{{Q + 4}}}}\left( {{B_R}} \right)}}
\label{3.7}%
\end{align}
and
\begin{align}
 &\quad {\left\| {{I_2}} \right\|_{{L^2}\left( {{B_R}} \right)}} \le 2{\left\| {\int_{{B_R}}
 {\left( {{A_0} - A} \right)\xi {D_0}u{D_0}{\Gamma _0}} d\zeta } \right\|_{{L^2}\left( {{B_R}} \right)}}
 + 2{\left\| {\int_{{B_R}} {{\Gamma _0}A{D_0}u{D_0}\xi } d\zeta } \right\|_{{L^2}\left( {{B_R}} \right)}} \nonumber\\
  & \le c{\left\| {{\Gamma _0}\left( {{D_0}\left( {\xi {D_0}u} \right)} \right)} \right\|_{{L^2}\left( {{B_R}} \right)}}
  + c{\left\| {{\Gamma _0}\left( {{D_0}u{D_0}\xi } \right)} \right\|_{{L^2}\left( {{B_R}} \right)}}\nonumber\\
  & \le c{\left\| {\xi {D_0}u} \right\|_{{L^{\frac{{2(Q + 2)}}{{Q + 4}}}}\left( {{B_R}} \right)}}
  + c{\left\| {{\Gamma _0}\left( {{D_0}u{D_0}\xi } \right)} \right\|_{{L^{\frac{{2(Q + 2)}}{Q}}}\left( {{B_R}} \right)}}
  {\left| {{B_R}} \right|^{\frac{1}{{Q + 2}}}} \nonumber\\
  & \le c{\left\| {\xi {D_0}u} \right\|_{{L^{\frac{{2(Q + 2)}}{{Q + 4}}}}\left( {{B_R}} \right)}}
  + c{\left\| {{D_0}u{D_0}\xi } \right\|_{{L^{\frac{{2(Q + 2)}}{{Q + 4}}}}\left( {{B_R}} \right)}}R \nonumber\\
& \le \frac{c}{{R - \rho }}{\left\| {{D_0}u}
\right\|_{{L^{\frac{{2(Q + 2)}}{{Q + 4}}}}\left( {{B_R}} \right)}}.
\label{3.8}%
\end{align}
Since $u$ is a weak solution to (\ref{eq1}), we infer that
$${I_3}(z) = \int_{{B_R}} {\left[ {f{D_0}\left( {\xi {\Gamma _0}}
\right) - g\xi {\Gamma _0}} \right]} d\zeta  = \int_{{B_R}} {\left[
{f\xi {D_0}{\Gamma _0} + f{\Gamma _0}{D_0}\xi  - g\xi {\Gamma _0}}
\right]} d\zeta $$
\noindent and
\begin{align}
& \quad {\left\| {{I_3}} \right\|_{{L^2}\left( {{B_R}} \right)}}\nonumber\\
& \le c{\left\| {\int_{{B_R}} {f\xi {D_0}{\Gamma _0}} d\zeta } \right\|_{{L^2}\left( {{B_R}} \right)}}
+ c{\left\| {\int_{{B_R}} {f{\Gamma _0}{D_0}\xi } d\zeta } \right\|_{{L^2}\left( {{B_R}} \right)}}
+ c{\left\| {\int_{{B_R}} {g\xi {\Gamma _0}} d\zeta } \right\|_{{L^2}\left( {{B_R}} \right)}} \nonumber\\
  & \le c{\left\| {{\Gamma _0}\left( {{D_0}\left( {f\xi } \right)} \right)} \right\|_{{L^2}\left( {{B_R}} \right)}}
  + c{\left\| {{\Gamma _0}\left( {f{D_0}\xi } \right)} \right\|_{{L^2}\left( {{B_R}} \right)}}
  + c{\left\| {{\Gamma _0}\left( {g\xi } \right)} \right\|_{{L^2}\left( {{B_R}} \right)}} \nonumber\\
 & \le c{\left\| {f\xi } \right\|_{{L^{\frac{{2(Q + 2)}}{{Q + 4}}}}\left( {{B_R}} \right)}}
  + cR\left( {{{\left\| {{\Gamma _0}\left( {f{D_0}\xi } \right)} \right\|}_{{L^{\frac{{2(Q + 2)}}{Q}}}\left( {{B_R}} \right)}}
  + c{{\left\| {{\Gamma _0}\left( {g\xi } \right)} \right\|}_{{L^{\frac{{2(Q + 2)}}{Q}}}\left( {{B_R}} \right)}}} \right) \nonumber\\
 & \le c{\left\| {f\xi } \right\|_{{L^{\frac{{2(Q + 2)}}{{Q + 4}}}}\left( {{B_R}} \right)}}
  + cR\left( {{{\left\| {f{D_0}\xi } \right\|}_{{L^{\frac{{2(Q + 2)}}{{Q + 4}}}}\left( {{B_R}} \right)}}
  + {{\left\| {g\xi } \right\|}_{{L^{\frac{{2(Q + 2)}}{{Q + 4}}}}\left( {{B_R}} \right)}}} \right) \nonumber\\
  & \le c{\left\| f \right\|_{{L^{\frac{{2(Q + 2)}}{{Q + 4}}}}\left( {{B_R}} \right)}}
  + cR\left( {\frac{c}{{R - \rho }}{{\left\| f \right\|}_{{L^{\frac{{2(Q + 2)}}{{Q + 4}}}}\left( {{B_R}} \right)}}
   + {{\left\| g \right\|}_{{L^{\frac{{2(Q + 2)}}{{Q + 4}}}}\left( {{B_R}} \right)}}} \right) \nonumber\\
& \le \frac{c}{{R - \rho }}\left( {{{\left\| f
\right\|}_{{L^{\frac{{2(Q + 2)}}{{Q + 4}}}}\left( {{B_R}} \right)}}
+ {{\left\| g \right\|}_{{L^{\frac{{2(Q + 2)}}{{Q + 4}}}}\left(
{{B_R}} \right)}}} \right).
\label{3.9}%
\end{align}
Inserting (3.7), (3.8) and (3.9) into (\ref{eq17}), it obtains (3.5).

\textbf{Theorem 3.3 }(Poincar\'{e} type inequality) Let $u \in
W_2^{1,1}\left( \Omega  \right)$ be a weak solution to (\ref{eq1}).
Then for any ${B_R} \subset \Omega$, $\rho  < R$, one has
\begin{equation}
\label{eq21} \int_{{B_\rho }} {{{\left| u \right|}^2}} dz \le
\frac{{c{R^4}}}{{{{\left( {R - \rho } \right)}^2}}}\int_{{B_R}}
{{{\left| {{D_0}u} \right|}^2}} dz + c{R^2}\int_{{B_R}} {\left(
{{{\left| g \right|}^2} + {{\left| f \right|}^2}} \right)} dz.
\end{equation}

\textbf{Proof:} Introduce two cutoff functions $\varsigma (x), \eta
(t) \in C_0^\infty \left( {{Q_R}} \right)$ satisfying
$$\varsigma (x) = 1(\left| x \right| < \rho), \quad \varsigma (x) = 0(\left| x \right| \ge
R),$$
$$0 \le \varsigma  \le 1, \quad \left| {{\partial _{{x_j}}}\varsigma } \right| \le \frac{c}{{R - \rho }}(j = 1, \ldots
,N);$$
\[\eta (t) = \left\{ {\begin{array}{*{20}{c}}
   {\frac{{2t - 2\left( {{t_0} - {R^2}/2} \right)}}{{{R^2} - {\rho ^2}}},t \in \left[ {{t_0} - \frac{{{R^2}}}{2},{t_0} - \frac{{{\rho ^2}}}{2}} \right),}  \\
   {{\kern 1pt} {\kern 1pt} {\kern 1pt} {\kern 1pt} {\kern 1pt} {\kern 1pt} {\kern 1pt} {\kern 1pt} {\kern 1pt} {\kern 1pt} {\kern 1pt} {\kern 1pt} {\kern 1pt} 1,{\kern 1pt} {\kern 1pt} {\kern 1pt} {\kern 1pt} {\kern 1pt} {\kern 1pt} {\kern 1pt} {\kern 1pt} {\kern 1pt} {\kern 1pt} {\kern 1pt} {\kern 1pt} {\kern 1pt} {\kern 1pt} {\kern 1pt} {\kern 1pt} {\kern 1pt} {\kern 1pt} {\kern 1pt} {\kern 1pt} {\kern 1pt} {\kern 1pt} {\kern 1pt} {\kern 1pt} {\kern 1pt} {\kern 1pt} {\kern 1pt} {\kern 1pt} {\kern 1pt} {\kern 1pt} {\kern 1pt} {\kern 1pt} {\kern 1pt} {\kern 1pt} {\kern 1pt} {\kern 1pt} {\kern 1pt} {\kern 1pt} {\kern 1pt} {\kern 1pt} {\kern 1pt} {\kern 1pt} {\kern 1pt} {\kern 1pt} {\kern 1pt} {\kern 1pt} t \in \left[ {{t_0} - \frac{{{\rho ^2}}}{2},{t_0} + \frac{{{R^2}}}{2}} \right].}  \\
\end{array}} \right.\]

Multiplying both sides of (\ref{eq1}) by $u{\varsigma ^2}(x)\eta
(t)$ and integrating on ${Q_R}^\prime  = {K_R} \times {S_R} \times
{I_R}^\prime$ (${I_R}^\prime  = \left[ {{t_0} - \frac{{{R^2}}}{2},s}
\right], s \le {t_0} + \frac{{{R^2}}}{2})$, we have
\begin{equation}
\label{eq22} \int_{{Q_R}^\prime } {\left[ { - A{D_0}u{D_0}\left(
{u{\varsigma ^2}\eta } \right) + xBu{\varsigma ^2}\eta Du -
u{\varsigma ^2}\eta {\partial _t}u} \right]} dz = \int_{{Q_R}^\prime
} {\left[ {gu{\varsigma ^2}\eta  - f{D_0}\left( {u{\varsigma ^2}\eta
} \right)} \right]} dz.
\end{equation}
Noting
\begin{equation}
\label{eq23} \int_{{Q_R}^\prime } {u{\varsigma ^2}\eta {\partial
_t}u} dz = \frac{1}{2}\int_{{Q_R}^\prime } {{\varsigma ^2}{{\left(
{{u^2}\eta } \right)}_t}} dz - \frac{1}{2}\int_{{Q_R}^\prime }
{{u^2}{\varsigma ^2}{\eta _t}} dz,
\end{equation}
\begin{equation}
\label{eq24}  \int_{{Q_R}^\prime } {xBu{\varsigma ^2}\eta Du} dz =
\frac{1}{2}\int_{{Q_R}^\prime } {xBD\left( {{u^2}{\varsigma ^2}\eta
} \right)} dz - \int_{{Q_R}^\prime } {xB{u^2}\varsigma \eta
D\varsigma } dz,
\end{equation}
it implies by inserting (\ref{eq23}) and (\ref{eq24}) into (\ref{eq22})
that
\begin{align}
& \quad \frac{1}{2}\int_{{Q_R}^\prime } {{u^2}{\varsigma ^2}{\eta _t}} dz \nonumber\\
  &= \int_{{Q_R}^\prime } {A{\varsigma ^2}\eta {D_0}u{D_0}u} dz{\rm{ + }}2\int_{{Q_R}^\prime }
  {Au\varsigma \eta {D_0}u{D_0}\varsigma } dz - \frac{1}{2}\int_{{Q_R}^\prime } {xBD\left( {{u^2}{\varsigma ^2}\eta } \right)} dz \nonumber\\
 &\quad + \int_{{Q_R}^\prime } {xB{u^2}\varsigma \eta D\varsigma } dz+ \frac{1}{2}\int_{{Q_R}^\prime } {{\varsigma ^2}{{\left( {{u^2}\eta } \right)}_t}} dz
 + \int_{{Q_R}^\prime } {gu{\varsigma ^2}\eta } dz \nonumber\\
 & \quad - \int_{{Q_R}^\prime } {f{\varsigma ^2}\eta {D_0}u} dz
 - 2\int_{{Q_R}^\prime } {fu\varsigma \eta {D_0}\varsigma } dz \nonumber\\
  & = \int_{{Q_R}^\prime } {A{\varsigma ^2}\eta {D_0}u{D_0}u} dz{\rm{ + }}2\int_{{Q_R}^\prime }
  {Au\varsigma \eta {D_0}u{D_0}\varsigma } dz - \int_{{Q_R}^\prime } {Y\left( {\frac{1}{2}{u^2}{\varsigma ^2}\eta } \right)} dz \nonumber\\
 & \quad + \int_{{Q_R}^\prime } {xB{u^2}\varsigma \eta D\varsigma } dz + \int_{{Q_R}^\prime } {gu{\varsigma ^2}\eta } dz -
\int_{{Q_R}^\prime } {f{\varsigma ^2}\eta {D_0}u} dz -
2\int_{{Q_R}^\prime } {fu\varsigma \eta {D_0}\varsigma } dz.
\label{3.14}%
\end{align}
By the divergence theorem and the property of $\varsigma$, it follows
\[\int_{{Q_R}^\prime } {Y\left( {\frac{1}{2}{u^2}{\varsigma ^2}\eta } \right)} dz = 0.\]
Hence we have by Young's inequality that
\begin{align}
& \quad \frac{1}{2}\int_{{Q_R}^\prime } {{u^2}{\varsigma ^2}{\eta _t}} dz \nonumber\\
 & \le \Lambda \int_{{Q_R}^\prime } {{{\left| {{D_0}u} \right|}^2}{\varsigma ^2}\eta } dz
  + \varepsilon \int_{{Q_R}^\prime } {{{\left| u \right|}^2}{{\left| {{D_0}\varsigma } \right|}^2}\eta } dz
  + {c_\varepsilon }\int_{{Q_R}^\prime } {{{\left| {{D_0}u} \right|}^2}{\varsigma ^2}\eta } dz \nonumber\\
 & \quad + c\int_{{Q_R}^\prime } {{{\left| u \right|}^2}\left| {D\varsigma } \right|\varsigma \eta } dz + {c_\varepsilon }\int_{{Q_R}^\prime } {{{\left| g \right|}^2}{\varsigma ^2}\eta } dz
 + \varepsilon \int_{{Q_R}^\prime } {{{\left| u \right|}^2}{\varsigma ^2}\eta } dz
 + {c_\varepsilon }\int_{{Q_R}^\prime } {{{\left| f \right|}^2}{\varsigma ^2}\eta } dz \nonumber\\
 &\quad + \varepsilon \int_{{Q_R}^\prime } {{{\left| {{D_0}u} \right|}^2}{\varsigma ^2}\eta } dz+ {c_\varepsilon }\int_{{Q_R}^\prime } {{{\left| f \right|}^2}{\varsigma ^2}\eta } dz
 + \varepsilon \int_{{Q_R}^\prime } {{{\left| u \right|}^2}{{\left| {{D_0}\varsigma } \right|}^2}\eta } dz \nonumber\\
& \le \int_{{Q_R}^\prime } {{{\left| u \right|}^2}\left(
{2\varepsilon {{\left| {{D_0}\varsigma } \right|}^2}\eta  +
c{{\left| {D\varsigma } \right|}^2}\varsigma \eta + \varepsilon
{\varsigma ^2}\eta } \right)} dz + c\int_{{Q_R}^\prime } {{{\left|
{{D_0}u} \right|}^2}{\varsigma ^2}\eta } dz \nonumber\\
&\quad + {c_\varepsilon }\int_{{Q_R}^\prime } {\left( {{{\left| g
\right|}^2} + {{\left| f \right|}^2}} \right){\varsigma ^2}\eta }
dz.\label{3.15}%
\end{align}
In the light of properties of $\varsigma ,\eta $ and (3.15), it yields
\begin{align}
& \quad \int_{{Q_\rho }} {{{\left| u \right|}^2}} dz \le \int_{{Q_R}^\prime } {{{\left| u \right|}^2}{\varsigma ^2}} dz
\le c\left( {{R^2} - {\rho ^2}} \right)\int_{{Q_R}^\prime } {{{\left| u \right|}^2}{\varsigma ^2}{\eta _t}} dz \nonumber\\
 & \le \left( {{R^2} - {\rho ^2}} \right)\int_{{Q_R}^\prime } {{{\left| u \right|}^2}\left( {2\varepsilon {{\left| {{D_0}\varsigma } \right|}^2}\eta
  + c{{\left| {D\varsigma } \right|}^2}\varsigma \eta  + \varepsilon {\varsigma ^2}\eta } \right)} dz \nonumber\\
 &\quad + c\left( {{R^2} - {\rho ^2}} \right)\int_{{Q_R}^\prime } {{{\left| {{D_0}u} \right|}^2}{\varsigma ^2}\eta } dz + {c_\varepsilon }\left( {{R^2} - {\rho ^2}} \right)\int_{{Q_R}^\prime } {\left( {{{\left| g \right|}^2}
 + {{\left| f \right|}^2}} \right){\varsigma ^2}\eta } dz \nonumber\\
 & \le \int_{{Q_R}} {{{\left| u \right|}^2}\left( {\frac{{2\varepsilon \left( {{R^2} - {\rho ^2}} \right)\eta }}{{{{\left( {R - \rho } \right)}^2}}}
 + \frac{{c\left( {{R^2} - {\rho ^2}} \right)\varsigma \eta }}{{{{\left( {R - \rho } \right)}^2}}}
 + \varepsilon \left( {{R^2} - {\rho ^2}} \right){\varsigma ^2}\eta } \right)} dz \nonumber\\
 &\quad + \frac{{c{R^2}{{\left( {R - \rho } \right)}^2}}}{{{{\left( {R - \rho } \right)}^2}}}\int_{{Q_R}} {{{\left| {{D_0}u} \right|}^2}} dz
 + {c_\varepsilon }{R^2}\int_{{Q_R}} {\left( {{{\left| g \right|}^2} + {{\left| f \right|}^2}} \right)} dz \nonumber\\
&\le {\theta _1}\int_{{Q_R}} {{{\left| u \right|}^2}} dz +
\frac{{c{R^4}}}{{{{\left( {R - \rho } \right)}^2}}}\int_{{Q_R}}
{{{\left| {{D_0}u} \right|}^2}} dz + {c_\varepsilon
}{R^2}\int_{{Q_R}} {\left( {{{\left| g \right|}^2} + {{\left| f
\right|}^2}} \right)} dz,
\label{3.16}%
\end{align}
\noindent where ${\theta _1} = \frac{{2\varepsilon \left( {{R^2} -
{\rho ^2}} \right)\eta }}{{{{\left( {R - \rho } \right)}^2}}} +
\frac{{c\left( {{R^2} - {\rho ^2}} \right)\varsigma \eta
}}{{{{\left( {R - \rho } \right)}^2}}} + \varepsilon \left( {{R^2} -
{\rho ^2}} \right)\varsigma \eta $. Choosing $\varepsilon$ small
enough, it ensures $0 < {\theta _1} < 1$ and we have from Lemma 2.9 that
\begin{equation}
\label{eq28} \int_{{Q_\rho }} {{{\left| u \right|}^2}} dz \le
\frac{{c{R^4}}}{{{{\left( {R - \rho } \right)}^2}}}\int_{{Q_R}}
{{{\left| {{D_0}u} \right|}^2}} dz + c{R^2}\int_{{Q_R}} {\left(
{{{\left| g \right|}^2} + {{\left| f \right|}^2}} \right)} dz.
\end{equation}
Now (\ref{eq28}) and ${B_{\rho /{c_0}}} \subset {Q_\rho } \subset
{Q_R} \subset {B_{{c_0}R}}$ imply (\ref{eq21}).


\section{Proof of Theorem 1.1}

\label{4}

Let us first describe a known result.

\textbf{Lemma 4.1 }(reverse H\"{o}lder inequality, [11]) Let $\hat
g$ and $\hat f$ be nonnegative functions on $\Omega$ and satisfy
\[\hat g \in {L^{\hat q}}(\Omega ), \quad \hat f \in {L^r}(\Omega ), \quad 1 < \hat q < r.
\]
If there exist constants $b_2$ and ${\theta _2}$ with ${b_2} > 1$
such that for any ${B_{2R}} \subset \Omega$, the inequality holds
\[\begin{array}{l}
 \quad \frac{1}{{\left| {{B_R}} \right|}}\int_{{B_R}} {{{\hat g}^{\hat q}}} dz \\
 \le {b_2}\left[ {{{\left( {\frac{1}{{\left| {{B_{4R/3}}} \right|}}\int_{{B_{4R/3}}} {\hat g} dz} \right)}^{\hat q}} + \frac{1}{{\left| {{B_{4R/3}}} \right|}}\int_{{B_{4R/3}}} {{{\hat f}^{\hat q}}} dz} \right] + {\theta _2}\frac{1}{{\left| {{B_{4R/3}}} \right|}}\int_{{B_{4R/3}}} {{{\hat g}^{\hat q}}} dz, \\
 \end{array}\]
then there exist positive constants ${\theta _0} = {\theta _0}(\hat
q,\Omega )$ and $\varepsilon _0$ such that if ${\theta _2} < {\theta
_0}$, then for any $\hat p \in \left[ {\hat q,\hat q + {\varepsilon _0}}
\right)$, it follows $\hat g \in L_{loc}^{\hat p}\left( \Omega \right)$
and
\begin{equation}
\label{eq29}  {\left( {\frac{1}{{\left| {{B_R}}
\right|}}\int_{{B_R}} {{{\hat g}^{\hat p}}} dz}
\right)^{\frac{1}{{\hat p}}}} \le c\left[ {{{\left(
{\frac{1}{{\left| {{B_{2R}}} \right|}}\int_{{B_{2R}}} {{{\hat
g}^{\hat q}}} dz} \right)}^{\frac{1}{{\hat q}}}}{\rm{ + }}{{\left(
{\frac{1}{{\left| {{B_{2R}}} \right|}}\int_{{B_{2R}}} {{{\hat
f}^{\hat p}}} dz} \right)}^{\frac{1}{{\hat p}}}}} \right],
\end{equation}
\noindent where $c$ and $\varepsilon _0$ depend on ${b_2},\hat
q,{\theta _2}$ and $Q$.

\textbf{Theorem 4.2} Let $u \in W_2^{1,1}\left( \Omega  \right)$ be
a weak solution to (\ref{eq1}) in $\Omega$. Then for any $p \in
\left[ {2,2 + \frac{{2Q}}{{Q + 2}}{\varepsilon _0}} \right)$, we
have ${D_0}u \in L_{loc}^p(\Omega )$ and for any ${B_R}
\subset {B_{2R}} \subset \Omega $,
\begin{align}
&\quad {\left( {\frac{1}{{\left| {{B_R}} \right|}}\int_{{B_R}}
{{{\left| {{D_0}u} \right|}^p}} dz}
\right)^{\frac{1}{p}}}\nonumber\\
& \le c\left[ {{{\left({\frac{1}{{\left| {{B_{2R}}}
\right|}}\int_{{B_{2R}}} {{{\left| {{D_0}u} \right|}^2}} dz}
\right)}^{\frac{1}{2}}} + {{\left( {\frac{1}{{\left| {{B_{2R}}}
\right|}}\int_{{B_{2R}}} {{{\left( {{{\left| g \right|}^2} +
{{\left| f \right|}^2}} \right)}^{\frac{p}{2}}}} dz}
\right)}^{\frac{1}{p}}}} \right].
\label{4.2}%
\end{align}

\textbf{Proof: }By using H\"{o}lder's inequality, it implies
\begin{align}
&\quad \int_{{B_{11R/9}}} {{{\left| {{D_0}u} \right|}^{\frac{{2(Q + 2)}}{{Q + 4}}}}} dz \nonumber\\
 & \le {\left( {\int_{{B_{11R/9}}} {{{\left| {{D_0}u} \right|}^2}} dz} \right)^{\frac{1}{2}}}
 {\left( {\int_{{B_{11R/9}}} {{{\left| {{D_0}u} \right|}^{\frac{{2Q}}{{Q + 4}}}}} dz} \right)^{\frac{1}{2}}} \nonumber\\
& \le {\left( {\int_{{B_{11R/9}}} {{{\left| {{D_0}u} \right|}^2}}
dz} \right)^{\frac{1}{2}}}{\left| {{B_{11R/9}}} \right|^{\frac{1}{{Q
+ 4}}}}{\left( {\int_{{B_{11R/9}}} {{{\left| {{D_0}u}
\right|}^{\frac{{2Q}}{{Q + 2}}}}} dz} \right)^{\frac{{Q + 2}}{{2(Q +
4)}}}}.\label{4.3}%
\end{align}
Combining (3.5) and (4.3), we get
\begin{align}
& \quad \int_{{B_{10R/9}}} {{{\left| u \right|}^2}} dz \nonumber\\
&\le \frac{c}{{{R^2}}} {\left[ {{{\left( {\int_{{B_{11R/9}}}
{{{\left| u \right|}^{\frac{{2(Q + 2)}}{{Q + 4}}}}} dz} \right)}
^{\frac{{Q + 4}}{{2(Q + 2)}}}} + {{\left( {\int_{{B_{11R/9}}}
{{{\left| {{D_0}u} \right|}
^{\frac{{2(Q + 2)}}{{Q + 4}}}}} dz} \right)}^{\frac{{Q + 4}}{{2(Q + 2)}}}}} \right]^2} \nonumber\\
 &\quad {\rm{ + }}\frac{c}{{{R^2}}}{\left[ {{{\left( {\int_{{B_{11R/9}}} {{{\left| f \right|}^{\frac{{2(Q + 2)}}{{Q + 4}}}}} dz} \right)}
 ^{\frac{{Q + 4}}{{2(Q + 2)}}}} + {{\left( {\int_{{B_{11R/9}}} {{{\left| g \right|}^{\frac{{2(Q + 2)}}{{Q + 4}}}}} dz} \right)}
 ^{\frac{{Q + 4}}{{2(Q + 2)}}}}} \right]^2} \nonumber\\
 & \le \frac{c}{{{R^2}}}{\left[ {{{\left| {{B_{11R/9}}} \right|}^{\frac{1}{{Q + 2}}}}
 {{\left( {\int_{{B_{11R/9}}} {{{\left| u \right|}^2}} dz} \right)}^{\frac{1}{2}}}} \right]^2} \nonumber\\
 &\quad + {\left[ {{{\left( {\int_{{B_{11R/9}}} {{{\left| {{D_0}u} \right|}^2}} dz} \right)}^{\frac{{Q + 4}}{{4(Q + 2)}}}}
 {{\left| {{B_{11R/9}}} \right|}^{\frac{1}{{2(Q + 2)}}}}{{\left( {\int_{{B_{11R/9}}} {{{\left| {{D_0}u} \right|}^{\frac{{2Q}}{{Q + 2}}}}} dz} \right)}
 ^{\frac{1}{4}}}} \right]^2} \nonumber\\
 &\quad {\rm{ + }}\frac{c}{{{R^2}}}{\left[ {{{\left| {{B_{11R/9}}} \right|}^{\frac{1}{{Q + 2}}}}{{\left( {\int_{{B_{11R/9}}}
 {{{\left| f \right|}^2}} dz} \right)}^{\frac{1}{2}}} + {{\left| {{B_{11R/9}}} \right|}^{\frac{1}{{Q + 2}}}}
 {{\left( {\int_{{B_{11R/9}}} {{{\left| g \right|}^2}} dz} \right)}^{\frac{1}{2}}}} \right]^2} \nonumber\\
 & \le c\int_{{B_{11R/9}}} {{{\left| u \right|}^2}} dz + \frac{c}{R}{\left( {\int_{{B_{11R/9}}} {{{\left| {{D_0}u} \right|}^2}} dz} \right)
 ^{\frac{{Q + 4}}{{2(Q + 2)}}}}{\left( {\int_{{B_{11R/9}}} {{{\left| {{D_0}u} \right|}^{\frac{{2Q}}{{Q + 2}}}}} dzdz} \right)^{\frac{1}{2}}} \nonumber\\
&\quad {\rm{ + }}c\int_{{B_{11R/9}}} {\left( {{{\left| f \right|}^2}
+ {{\left| g \right|}^2}} \right)} dz.
\label{4.4}%
\end{align}
Noting (\ref{eq13}), (\ref{eq21}) and (4.4), it follows
\[\begin{array}{l}
 \quad \int_{{B_R}} {{{\left| {{D_0}u} \right|}^2}} dz \\
  \le \frac{c}{{{R^2}}}\int_{{B_{11R/9}}} {{{\left| u \right|}^2}} dz + \frac{c}{{{R^3}}}{\left( {\int_{{B_{11R/9}}} {{{\left| {{D_0}u} \right|}^2}} dz} \right)^{\frac{{Q + 4}}{{2(Q + 2)}}}}{\left( {\int_{{B_{11R/9}}} {{{\left| {{D_0}u} \right|}^{\frac{{2Q}}{{Q + 2}}}}} dz} \right)^{\frac{1}{2}}} \\
 {\kern 1pt} {\kern 1pt} {\kern 1pt} {\kern 1pt} {\kern 1pt} {\kern 1pt} {\kern 1pt} {\kern 1pt} {\kern 1pt} {\rm{ + }}\frac{c}{{{R^2}}}\int_{{B_{11R/9}}} {\left( {{{\left| f \right|}^2} + {{\left| g \right|}^2}} \right)} dz + c\int_{{B_{10R/9}}} {\left( {{{\left| g \right|}^2} + {{\left| f \right|}^2}} \right)} dz \\
 \le c\int_{{B_{4R/3}}} {{{\left| {{D_0}u} \right|}^2}} dz\\
 {\kern 1pt} {\kern 1pt} {\kern 1pt} {\kern 1pt} {\kern 1pt} {\kern 1pt} {\kern 1pt} {\kern 1pt}  + \frac{c}{{{R^3}}}{\left| {{B_{4R/3}}} \right|^{\frac{{Q + 3}}{{Q + 2}}}}{\left( {\frac{1}{{\left| {{B_{4R/3}}} \right|}}\int_{{B_{4R/3}}} {{{\left| {{D_0}u} \right|}^2}} dz} \right)^{\frac{{Q + 4}}{{2(Q + 2)}}}}{\left( {\frac{1}{{\left| {{B_{4R/3}}} \right|}}\int_{{B_{4R/3}}} {{{\left| {{D_0}u} \right|}^{\frac{{2Q}}{{Q + 2}}}}} dz} \right)^{\frac{1}{2}}}\\
\quad + \frac{c}{{{R^2}}}\int_{{B_{4R/3}}} {\left( {{{\left| g \right|}^2} + {{\left| f \right|}^2}} \right)} dz \\
 \end{array}\]
and hence
\begin{align}
& \quad \frac{1}{{\left| {{B_R}} \right|}}\int_{{B_R}} {{{\left| {{D_0}u} \right|}^2}} dz \nonumber\\
  & \le \frac{c}{{\left| {{B_{4R/3}}} \right|}}\int_{{B_{4R/3}}} {{{\left| {{D_0}u} \right|}^2}} dz \nonumber\\
 &\quad + \varepsilon \left( {\frac{1}{{\left| {{B_{4R/3}}} \right|}}\int_{{B_{4R/3}}} {{{\left| {{D_0}u} \right|}^2}} dz} \right)
 + {c_\varepsilon }{R^{ - \frac{{4(Q + 2)}}{Q}}}{\left( {\frac{1}{{\left| {{B_{4R/3}}} \right|}}
 \int_{{B_{4R/3}}} {{{\left| {{D_0}u} \right|}^{\frac{{2Q}}{{Q + 2}}}}} dz} \right)^{\frac{{Q + 2}}{Q}}} \nonumber\\
  & \quad + \frac{c}{{{R^2}}}\frac{1}{{\left| {{B_{4R/3}}} \right|}}\int_{{B_{4R/3}}} {\left( {{{\left| g \right|}^2}
  + {{\left| f \right|}^2}} \right)} dz \nonumber\\
  & \le c\left( {\frac{1}{{\left| {{B_{4R/3}}}
\right|}}\int_{{B_{4R/3}}} {{{\left| {{D_0}u} \right|}^2}} dz}
\right) + {c_\varepsilon }{R^{ - \frac{{4(Q + 2)}}{Q}}}{\left( {\frac{1}{{\left| {{B_{4R/3}}} \right|}}\int_{{B_{4R/3}}}
  {{{\left| {{D_0}u} \right|}^{\frac{{2Q}}{{Q + 2}}}}} dz} \right)^{\frac{{Q + 2}}{Q}}} \nonumber\\
&\quad+ \frac{c}{{{R^2}}}\frac{1}{{\left| {{B_{4R/3}}} \right|}}\int_{{B_{4R/3}}}
{\left( {{{\left| g \right|}^2} + {{\left| f \right|}^2}} \right)} dz.\label{4.5}%
\end{align}
Let $\hat g = {\left| {{D_0}u} \right|^{\tilde q}}$, $\tilde q =
\frac{{2Q}}{{Q + 2}}$, $\hat q = \frac{2}{{\tilde q}} = \frac{{Q +
2}}{Q} > 1$, $\hat f = {\left( {{{\left| g \right|}^2} + {{\left| f
\right|}^2}} \right)^{\frac{Q}{{Q + 2}}}}$, then we rewrite (4.5) in the
form
\begin{align}
& \quad \frac{1}{{\left| {{B_R}} \right|}}\int_{{B_R}} {{{\hat
g}^{\hat q}}} dz \nonumber\\
& \le c\left[ {{{\left( {\frac{1}{{\left|
{{B_{4R/3}}} \right|}}\int_{{B_{4R/3}}} {\hat g} dz} \right)}^{\hat
q}} + \frac{1}{{\left| {{B_{4R/3}}} \right|}}\int_{{B_{4R/3}}}
{{{\hat f}^{\hat q}}} dz} \right] + \frac{c}{{\left| {{B_{4R/3}}}
\right|}}\int_{{B_{4R/3}}} {{{\hat g}^{\hat q}}} dz.
\label{4.6}%
\end{align}
It shows from Lemma 4.1 that for any $\hat p \in \left[ {\hat
q,\hat q + {\varepsilon _0}} \right)$,
\[{\left( {\frac{1}{{\left| {{B_R}} \right|}}\int_{{B_R}} {{{\hat g}^{\hat p}}} dz} \right)^{1/\hat p}} \le c\left[ {{{\left( {\frac{1}{{\left| {{B_{2R}}} \right|}}\int_{{B_{2R}}} {{{\hat g}^{\hat q}}} dz} \right)}^{1/\hat q}} + {{\left( {\frac{1}{{\left| {{B_{2R}}} \right|}}\int_{{B_{2R}}} {{{\hat f}^{\hat p}}} dz} \right)}^{1/\hat p}}} \right],\]
which means
\begin{align}
&\quad {\left( {\frac{1}{{\left| {{B_R}} \right|}}\int_{{B_R}}
{{{\left| {{D_0}u} \right|}^{\tilde q}}^{\hat p}} dz}
\right)^{\frac{1}{{\hat p}}}}\nonumber\\
 &\le c\left[ {{{\left(
{\frac{1}{{\left| {{B_{2R}}} \right|}}\int_{{B_{2R}}} {{{\left|
{{D_0}u} \right|}^2}} dz} \right)}^{\frac{Q}{{Q + 2}}}} + {{\left(
{\frac{1}{{\left| {{B_{2R}}} \right|}}\int_{{B_{2R}}} {{{\left(
{{{\left| g \right|}^2} + {{\left| f \right|}^2}}
\right)}^{\frac{{\hat p\tilde q}}{2}}}} dz} \right)}^{\frac{1}{{\hat
p}}}}} \right].\label{4.7}%
\end{align}
Setting $p = \hat p\tilde q \in \left[ {2,2 + \frac{{2Q}}{{Q +
2}}{\varepsilon _0}} \right)$, we finish the proof.

\textbf{Proof Theorem 1.1: } The conclusion follows from Theorem 4.2 and the cutoff function
technique.


\section{Homogeneous ultraparabolic equation}

\label{5}

In this section, we consider the following homogeneous
ultraparabolic equation
\begin{equation}
\label{eq36} div\left( {A{D_0}u} \right) + Yu = 0.
\end{equation}
To obtain $L^p$ estimates for gradients of weak solutions to (5.1), we divide (\ref{eq36}) into two parts. Let $v$ be a
weak solution to the following Dirichlet boundary value condition to the
homogeneous ultraparabolic equation with constant principal part:
\begin{equation}
\label{eq37} \left\{ \begin{array}{l}
 div\left( {{A_R}{D_0}v} \right) + Yv = 0,   in{\kern 1pt} {\kern 1pt} {\kern 1pt} {B_R}, \\
\quad \quad \quad  v = u, \quad  \quad  \quad  on \quad  {\partial _p}{B_R}. \\
 \end{array} \right.
\end{equation}
Then $w = u - v$ satisfies the Dirichlet
boundary value condition to the nonhomogeneous ultraparabolic equation
with constant principal part:
\begin{equation}
\label{eq38} \left\{ \begin{array}{l}
 div\left( {{A_R}{D_0}w} \right) + Yw = div\left( {\left( {{A_R} - A} \right){D_0}u} \right),{\kern 1pt} {\kern 1pt} {\kern 1pt} {\kern 1pt} {\kern 1pt} {\kern 1pt} {\kern 1pt} {\kern 1pt} in{\kern 1pt} {\kern 1pt} {\kern 1pt} {B_R}, \\
 {\kern 1pt} {\kern 1pt} {\kern 1pt} {\kern 1pt} {\kern 1pt} {\kern 1pt} {\kern 1pt} {\kern 1pt} {\kern 1pt} {\kern 1pt} {\kern 1pt} {\kern 1pt} {\kern 1pt} {\kern 1pt} {\kern 1pt} {\kern 1pt} {\kern 1pt} {\kern 1pt} {\kern 1pt} {\kern 1pt} {\kern 1pt} {\kern 1pt} {\kern 1pt} {\kern 1pt} {\kern 1pt} {\kern 1pt} {\kern 1pt} {\kern 1pt} {\kern 1pt} {\kern 1pt} {\kern 1pt} {\kern 1pt} {\kern 1pt} {\kern 1pt} {\kern 1pt} {\kern 1pt} {\kern 1pt} {\kern 1pt} {\kern 1pt} {\kern 1pt} {\kern 1pt} {\kern 1pt} {\kern 1pt} {\kern 1pt} {\kern 1pt} {\kern 1pt} {\kern 1pt} {\kern 1pt} {\kern 1pt} {\kern 1pt} {\kern 1pt} {\kern 1pt} {\kern 1pt} {\kern 1pt} {\kern 1pt} {\kern 1pt} {\kern 1pt} {\kern 1pt} {\kern 1pt} {\kern 1pt} {\kern 1pt} {\kern 1pt} {\kern 1pt} {\kern 1pt} {\kern 1pt} {\kern 1pt} {\kern 1pt} {\kern 1pt} {\kern 1pt} {\kern 1pt} {\kern 1pt} {\kern 1pt} {\kern 1pt} {\kern 1pt} w = 0,{\kern 1pt} {\kern 1pt} {\kern 1pt} {\kern 1pt} {\kern 1pt} {\kern 1pt} {\kern 1pt} {\kern 1pt} {\kern 1pt} {\kern 1pt} {\kern 1pt} {\kern 1pt} {\kern 1pt} {\kern 1pt} {\kern 1pt} {\kern 1pt} {\kern 1pt} {\kern 1pt} {\kern 1pt} {\kern 1pt} {\kern 1pt} {\kern 1pt} {\kern 1pt} {\kern 1pt} {\kern 1pt} {\kern 1pt} {\kern 1pt} {\kern 1pt} {\kern 1pt} {\kern 1pt} {\kern 1pt} {\kern 1pt} {\kern 1pt} {\kern 1pt} {\kern 1pt} {\kern 1pt} {\kern 1pt} {\kern 1pt} {\kern 1pt} {\kern 1pt} {\kern 1pt} {\kern 1pt} {\kern 1pt} {\kern 1pt} {\kern 1pt} {\kern 1pt} {\kern 1pt} {\kern 1pt} {\kern 1pt} {\kern 1pt} {\kern 1pt} {\kern 1pt} {\kern 1pt} {\kern 1pt} {\kern 1pt} {\kern 1pt} {\kern 1pt} {\kern 1pt} {\kern 1pt} {\kern 1pt} {\kern 1pt} {\kern 1pt} {\kern 1pt} {\kern 1pt} {\kern 1pt} {\kern 1pt} {\kern 1pt} {\kern 1pt} {\kern 1pt} {\kern 1pt} {\kern 1pt} {\kern 1pt} {\kern 1pt} {\kern 1pt} {\kern 1pt} {\kern 1pt} {\kern 1pt} {\kern 1pt} {\kern 1pt} {\kern 1pt} on{\kern 1pt} {\kern 1pt} {\kern 1pt} {\kern 1pt} \partial {B_R}{\rm{,}} \\
 \end{array} \right.
\end{equation}
\noindent where ${A_R} = \frac{1}{{\left| {{B_R}}
\right|}}\int_{{B_R}} A dz$.

\textbf{Lemma 5.1} Let $v \in W_2^{1,1}\left( \Omega  \right)$ be a
weak solution to (\ref{eq37}). Then for any ${B_R} \subset \Omega$,
one has
\begin{equation}
\label{eq39}  \mathop {\sup }\limits_{{B_{R{\rm{/}}2}}} {\left| v
\right|^2} \le \frac{c}{{{R^{Q{\rm{ + }}2}}}}\int_{{B_R}} {{{\left|
v \right|}^2}} dz{\rm{.}}
\end{equation}

\textbf{Proof}: It is true from Corollary 1.4 of [22].

Furthermore, we have

\textbf{Lemma 5.2} Let $v \in W_2^{1,1}\left( \Omega  \right)$ be a
weak solution to (\ref{eq37}). Then for any ${B_R} \subset \Omega$,
$\rho  < R$, it follows
\begin{equation}
\label{eq40}  \int_{{B_\rho }} {{{\left| v \right|}^2}} dz \le
c{\left( {\frac{\rho }{R}} \right)^{Q + 2}}\int_{{B_R}} {{{\left| v
\right|}^2}} dz.
\end{equation}

\textbf{Proof}: When $\frac{R}{2} \le \rho  < R$, the result is
evident. Now it is enough to treat the case $\rho  < \frac{R}{2}$.
But by Lemma 5.1, it yields
\[\begin{array}{l}
 \quad \int_{{B_\rho }} {{{\left| v \right|}^2}} dz
 \le \left| {{B_\rho }} \right|\mathop {\sup }\limits_{{B_\rho }} {\left| v \right|^2}
 \le \left| {{B_\rho }} \right|\mathop {\sup }\limits_{{B_{R{\rm{/}}2}}} {\left| v \right|^2} \\
  \le \left| {{B_\rho }} \right|\frac{c}{{{R^{Q{\rm{ + }}2}}}}\int_{{B_R}} {{{\left| v \right|}^2}} dz
  \le c{\left( {\frac{\rho }{R}} \right)^{Q{\rm{ + }}2}}\int_{{B_R}} {{{\left| v \right|}^2}} dz{\rm{.}} \\
 \end{array}\]

On the gradient of $v$, we have

\textbf{Lemma 5.3} Let $v \in W_2^{1,1}\left( \Omega  \right)$ be a
weak solution to (\ref{eq37}). Then for any ${B_R} \subset \Omega$,
$\rho  < R$, it follows
\begin{equation}
\label{eq41}  \int_{{B_\rho }} {{{\left| {{D_0}v} \right|}^2}} dz
\le c{\left( {\frac{\rho }{R}} \right)^Q}\int_{{B_R}} {{{\left|
{{D_0}v} \right|}^2}} dz.
\end{equation}

\textbf{Proof}: Combining Theorem 3.1, Theorem 3.3 ($g{\rm{ =
}}f{\rm{ = }}0$) and (\ref{eq40}), we arrive at
\[\begin{array}{l}
 \quad \int_{{B_{\rho {\rm{/}}2}}} {{{\left| {{D_0}v} \right|}^2}} dz \le \frac{c}{{{\rho ^2}}}\int_{{B_\rho }} {{{\left| v \right|}^2}} dz \le \frac{c}{{{\rho ^2}}}{\left( {\frac{\rho }{R}} \right)^{Q{\rm{ + }}2}}\int_{{B_R}} {{{\left| v \right|}^2}} dz \\
  \le \frac{c}{{{\rho ^2}}}{\left( {\frac{\rho }{R}} \right)^{Q{\rm{ + }}2}}{R^2}\int_{{B_{2R}}} {{{\left| {{D_0}v} \right|}^2}} dz \le c{\left( {\frac{\rho }{R}} \right)^Q}\int_{{B_{2R}}} {{{\left| {{D_0}v} \right|}^2}} dz{\rm{.}} \\
 \end{array}\]

\textbf{Lemma 5.4} Let $v \in W_2^{1,1}\left( \Omega  \right)$ be a
weak solution to (\ref{eq37}). Then for any $p \in \left[ {2,2 +
\frac{{2Q}}{{Q + 2}}{\varepsilon _0}} \right)$, ${B_R} \subset
\Omega$, $\rho  < R$, we have
\begin{equation}
\label{eq42} \int_{{B_\rho }} {{{\left| {{D_0}v} \right|}^p}} dz \le
c{\left( {\frac{\rho }{R}} \right)^{Q + 2 - p}}\int_{{B_R}}
{{{\left| {{D_0}v} \right|}^p}} dz.
\end{equation}

\textbf{Proof}: By Theorem 4.2($g{\rm{ = }}f{\rm{ = }}0$) and
(\ref{eq41}),
\[{\left( {\frac{1}{{\left| {{B_{\rho /2}}} \right|}}\int_{{B_{\rho {\rm{/}}2}}} {{{\left| {{D_0}v} \right|}^p}} dz} \right)^{\frac{1}{p}}} \le c{\left( {\frac{1}{{\left| {{B_\rho }} \right|}}\int_{{B_\rho }} {{{\left| {{D_0}v} \right|}^2}} dz} \right)^{\frac{1}{2}}} \le c{\left( {\frac{1}{{\left| {{B_\rho }} \right|}}{{\left( {\frac{\rho }{R}} \right)}^Q}\int_{{B_R}} {{{\left| {{D_0}v} \right|}^2}} dz} \right)^{\frac{1}{2}}}{\rm{.}}\]
From H\"{o}lder's inequality, it implies
\[\begin{array}{l}
 \quad \int_{{B_{\rho {\rm{/}}2}}} {{{\left| {{D_0}v} \right|}^p}} dz \le c\left| {{B_{\rho /2}}} \right|{\left( {\frac{1}{{\left| {{B_\rho }} \right|}}{{\left( {\frac{\rho }{R}} \right)}^Q}\int_{{B_R}} {{{\left| {{D_0}v} \right|}^2}} dz} \right)^{\frac{p}{2}}} \\
  \le c\left| {{B_{\frac{\rho }{2}}}} \right|\frac{1}{{{{\left| {{B_\rho }} \right|}^{\frac{p}{2}}}}}{\left( {\frac{\rho }{R}} \right)^{\frac{{pQ}}{2}}}{\left| {{B_R}} \right|^{\frac{{p - 2}}{2}}}\int_{{B_R}} {{{\left| {{D_0}v} \right|}^p}} dz \le c{\left( {\frac{{\left| {{B_\rho }} \right|}}{{\left| {{B_R}} \right|}}} \right)^{\frac{{2 - p}}{2}}}{\left( {\frac{\rho }{R}} \right)^{\frac{{pQ}}{2}}}\int_{{B_R}} {{{\left| {{D_0}v} \right|}^p}} dz \\
  \le c{\left( {\frac{\rho }{R}} \right)^{Q + 2 - p}}\int_{{B_R}} {{{\left| {{D_0}v} \right|}^p}} dz \\
 \end{array}\]
and the proof is ended.

The main result of this section is

\textbf{Theorem 5.5} Let $u \in W_2^{1,1}\left( \Omega  \right)$ be
a weak solution to (\ref{eq36}). Then for any $p \in \left[ {2,2 +
\frac{{2Q}}{{Q + 2}}{\varepsilon _0}} \right)$, $\varepsilon _0$ is the constant in Theorem 1.1, $\frac{{p - 2}}{p}(Q
+ 2) < \mu  < Q$, ${B_R} \subset \Omega$, $\rho  < R$, one has

\begin{equation}
\label{eq43} \int_{{B_\rho }} {{{\left| {{D_0}u} \right|}^p}} dz \le
c{\left( {\frac{\rho }{R}} \right)^{\frac{{2(Q + 2) - p(Q + 2 - \mu
)}}{2}}}\int_{{B_R}} {{{\left| {{D_0}u} \right|}^p}} dz.
\end{equation}

\textbf{Proof:} When $\frac{1}{2}R \le \rho  < R$, (5.8) is
clearly true. The remainder is to treat the case $\rho  < \frac{1}{2}R$.

Multiplying both sides of (\ref{eq38}) by $w$ and integrating on
$B_R$, it observes
\begin{equation}
\label{eq44}  - \int_{{B_R}} {{A_R}{D_0}w{D_0}w} dz{\rm{ +
}}\int_{{B_R}} {wYw} dz =  - \int_{{B_R}} {\left( {{A_R} - A}
\right){D_0}u{D_0}w} dz,
\end{equation}
and from the divergence theorem,
$$\int_{{B_R}} {wYw} dz = \frac{1}{2}\int_{{B_R}} {Y\left( {{w^2}} \right)} dz = 0.$$
By (H1) and Young's inequality, we have from (5.9) that
\begin{equation}
\label{eq45} {\Lambda ^{ - 1}}\int_{{B_R}} {{{\left| {{D_0}w}
\right|}^2}} dz \le {c_\varepsilon }\int_{{B_R}} {{{\left| {{A_R} -
A} \right|}^2}{{\left| {{D_0}u} \right|}^2}} dz + \varepsilon
\int_{{B_R}} {{{\left| {{D_0}w} \right|}^2}} dz.
\end{equation}
Choosing $\varepsilon$ small enough such that ${\Lambda ^{ - 1}} -
\varepsilon  > 0$, then (\ref{eq45}) implies
\begin{align}
&\quad  \int_{{B_R}} {{{\left| {{D_0}w} \right|}^2}} dz \le c\int_{{B_R}} {{{\left| {{A_R} - A} \right|}^2}{{\left| {{D_0}u} \right|}^2}} dz \nonumber\\
  &\le c{\left( {\int_{{B_R}} {{{\left| {{A_R} - A} \right|}^{\frac{{2p}}{{p - 2}}}}} dz} \right)^{\frac{{p - 2}}{p}}}
  {\left( {\int_{{B_R}} {{{\left| {{D_0}u} \right|}^p}} dz} \right)^{\frac{2}{p}}} \nonumber\\
& \le c{\left( {\left| {{B_R}} \right|{\eta _R}\left( {{a_{ij}}}
\right)} \right)^{\frac{{p - 2}}{p}}}{\left( {\int_{{B_R}} {{{\left|
{{D_0}u} \right|}^p}} dz} \right)^{\frac{2}{p}}}{\rm{}}
\label{5.11}%
\end{align}
and applying (\ref{eq41}) and (5.11) leads to
\begin{align}
&\quad \int_{{B_{2\rho }}} {{{\left| {{D_0}u} \right|}^2}} dz \le 2\int_{{B_{2\rho }}} {{{\left| {{D_0}v} \right|}^2}} dz
+ 2\int_{{B_{2\rho }}} {{{\left| {{D_0}w} \right|}^2}} dz \nonumber\\
  & \le c{\left( {\frac{\rho }{R}} \right)^Q}\int_{{B_R}} {{{\left| {{D_0}v} \right|}^2}} dz + c\int_{{B_R}} {{{\left| {{D_0}w} \right|}^2}} dz \nonumber\\
  &\le c{\left( {\frac{\rho }{R}} \right)^Q}\int_{{B_R}} {{{\left| {{D_0}u} \right|}^2}} dz + c\int_{{B_R}} {{{\left| {{D_0}w} \right|}^2}} dz \nonumber\\
 & \le c{\left( {\frac{\rho }{R}} \right)^Q}{\left| {{B_R}} \right|^{\frac{{p - 2}}{p}}}
 {\left( {\int_{{B_R}} {{{\left| {{D_0}u} \right|}^p}} dz} \right)^{\frac{2}{p}}}{\rm{ + }}
 c{\left( {\left| {{B_R}} \right|{\eta _R}\left( {{a_{ij}}} \right)} \right)^{\frac{{p - 2}}{p}}}
 {\left( {\int_{{B_R}} {{{\left| {{D_0}u} \right|}^p}} dz} \right)^{\frac{2}{p}}} \nonumber\\
& \le c\left[ {{{\left( {\frac{\rho }{R}} \right)}^Q}{\rm{ +
}}{{\left( {{\eta _R}\left( {{a_{ij}}} \right)} \right)}^{\frac{{p -
2}}{p}}}} \right]{\left( {{{\left| {{B_R}} \right|}^{\frac{{p -
2}}{2}}}\int_{{B_R}} {{{\left| {{D_0}u} \right|}^p}} dz}
\right)^{\frac{2}{p}}}{\rm{.}}
\label{5.12}%
\end{align}
It shows owing to Theorem 4.2 ($g = f = 0)$ that
\begin{align}
&\quad  {\left( {{{\left| {{B_\rho }} \right|}^{\frac{{p -
2}}{2}}}\int_{{B_\rho }} {{{\left| {{D_0}u} \right|}^p}} dz}
\right)^{\frac{2}{p}}} \le c\int_{{B_{2\rho }}} {{{\left| {{D_0}u}
\right|}^2}} dz \nonumber\\
&\le c\left[ {{{\left( {\frac{\rho }{R}}
\right)}^Q}{\rm{ + }}{{\left( {{\eta _R}\left( {{a_{ij}}} \right)}
\right)}^{\frac{{p - 2}}{p}}}} \right]{\left( {{{\left| {{B_R}}
\right|}^{\frac{{p - 2}}{2}}}\int_{{B_R}} {{{\left| {{D_0}u}
\right|}^p}} dz} \right)^{\frac{2}{p}}}{\rm{.}}
\label{5.13}%
\end{align}
Denoting $H(\rho ) = {\left( {{{\left| {{B_\rho }}
\right|}^{\frac{{p - 2}}{2}}}\int_{{B_\rho }} {{{\left| {{D_0}u}
\right|}^p}} dz} \right)^{\frac{2}{p}}}$, $H(R) = {\left( {{{\left|
{{B_R}} \right|}^{\frac{{p - 2}}{2}}}\int_{{B_R}} {{{\left| {{D_0}u}
\right|}^p}} dz} \right)^{\frac{2}{p}}}$, ${a_1} = Q$, ${B_1} = 0$ in Lemma 2.10,
we know that there exists ${b_1} = \mu (\frac{{p -
2}}{p}\left( {Q + 2} \right) < \mu  < Q)$ such that
\begin{equation}
\label{eq49} {\left( {{{\left| {{B_\rho }} \right|}^{\frac{{p -
2}}{2}}}\int_{{B_\rho }} {{{\left| {{D_0}u} \right|}^p}} dz}
\right)^{\frac{2}{p}}} \le c{\left( {\frac{\rho }{R}} \right)^\mu
}{\left( {{{\left| {{B_R}} \right|}^{\frac{{p - 2}}{2}}}\int_{{B_R}}
{{{\left| {{D_0}u} \right|}^p}} dz} \right)^{\frac{2}{p}}}.
\end{equation}
Inserting $\frac{{\left| {{B_R}} \right|}}{{\left| {{B_\rho }}
\right|}} \le c{\left( {\frac{\rho }{R}} \right)^{ - Q - 2}}$ into
(\ref{eq49}), it attains (\ref{eq43}).


\section{Proof of Theorem 1.2}

\label{6}
Based on the discussion in the preceding section, let $v$ be a weak solution to the
following problem
\begin{equation}
\label{eq50} \left\{ \begin{array}{l}
 div\left( {A{D_0}v} \right) + Yv = 0,\quad  in \quad {B_R}, \\
\quad \quad   v{\rm{ = }}u, \quad  on \quad \partial {B_R}, \\
 \end{array} \right.
\end{equation}
\noindent then $w=u-v$ satisfies
\begin{equation}
\label{eq51} \left\{ \begin{array}{l}
 div\left( {A{D_0}w} \right) + Yw = g + divf,{\kern 1pt} {\kern 1pt} {\kern 1pt} {\kern 1pt} {\kern 1pt} {\kern 1pt} {\kern 1pt} {\kern 1pt} {\kern 1pt} in{\kern 1pt} {\kern 1pt} {\kern 1pt} {\kern 1pt} {B_R}, \\
 {\kern 1pt} {\kern 1pt} {\kern 1pt} {\kern 1pt} {\kern 1pt} {\kern 1pt} {\kern 1pt} {\kern 1pt} {\kern 1pt} {\kern 1pt} {\kern 1pt} {\kern 1pt} {\kern 1pt} {\kern 1pt} {\kern 1pt} {\kern 1pt} {\kern 1pt} {\kern 1pt} {\kern 1pt} {\kern 1pt} {\kern 1pt} {\kern 1pt} {\kern 1pt} {\kern 1pt} {\kern 1pt} {\kern 1pt} {\kern 1pt} {\kern 1pt} {\kern 1pt} {\kern 1pt} {\kern 1pt} {\kern 1pt} {\kern 1pt} {\kern 1pt} {\kern 1pt} {\kern 1pt} {\kern 1pt} {\kern 1pt} {\kern 1pt} {\kern 1pt} {\kern 1pt} {\kern 1pt} {\kern 1pt} {\kern 1pt} {\kern 1pt} {\kern 1pt} {\kern 1pt} {\kern 1pt} {\kern 1pt} {\kern 1pt} {\kern 1pt} {\kern 1pt} {\kern 1pt} {\kern 1pt} {\kern 1pt} {\kern 1pt} {\kern 1pt} {\kern 1pt} {\kern 1pt} {\kern 1pt} {\kern 1pt} {\kern 1pt} {\kern 1pt} {\kern 1pt} {\kern 1pt} {\kern 1pt} {\kern 1pt} {\kern 1pt} {\kern 1pt} {\kern 1pt} w = 0,{\kern 1pt} {\kern 1pt} {\kern 1pt} {\kern 1pt} {\kern 1pt} {\kern 1pt} {\kern 1pt} {\kern 1pt} {\kern 1pt} {\kern 1pt} {\kern 1pt} {\kern 1pt} {\kern 1pt} {\kern 1pt} {\kern 1pt} {\kern 1pt} {\kern 1pt} {\kern 1pt} {\kern 1pt} {\kern 1pt} {\kern 1pt} {\kern 1pt} {\kern 1pt} {\kern 1pt} {\kern 1pt} {\kern 1pt} {\kern 1pt} {\kern 1pt} {\kern 1pt} {\kern 1pt} {\kern 1pt} {\kern 1pt} {\kern 1pt} {\kern 1pt} {\kern 1pt} on{\kern 1pt} {\kern 1pt} {\kern 1pt} {\kern 1pt} \partial {B_R}. \\
 \end{array} \right.
\end{equation}

\textbf{Theorem 6.1 } Let $w \in W_{2,0}^{1,1}\left( \Omega \right)$
be a weak solution to (\ref{eq51}). Then for any ${B_{2R}} \subset
\Omega$, one has
\begin{equation}
\label{eq52} \int_{{B_R}} {{{\left| {{D_0}w} \right|}^2}} dz \le
c\int_{{B_{2R}}} {\left( {{{\left| g \right|}^2} + {{\left| f
\right|}^2}} \right)} dz.
\end{equation}

\textbf{Proof:} Multiplying both sides of (\ref{eq51}) by $w$ and
integrating on $B_R$,
\begin{equation}
\label{eq53}  - \int_{{B_R}} {A{D_0}w{D_0}w} dz{\rm{ +
}}\int_{{B_R}} {wYw} dz = \int_{{B_R}} {gw} dz - \int_{{B_R}}
{f{D_0}w} dz.
\end{equation}
By (H1), the divergence theorem and Young's inequality, we have
\begin{equation}
\label{eq54} {\Lambda ^{ - 1}}\int_{{B_R}} {{{\left| {{D_0}w}
\right|}^2}} dz \le {c_\varepsilon }\int_{{B_R}} {{{\left| g
\right|}^2}} dz + \varepsilon \int_{{B_R}} {{{\left| w \right|}^2}}
dz + {c_\varepsilon }\int_{{B_R}} {{{\left| f \right|}^2}} dz +
\varepsilon \int_{{B_R}} {{{\left| {{D_0}w} \right|}^2}} dz.
\end{equation}
Since by using (\ref{eq21}),
\begin{equation}
\label{eq55} \int_{{B_R}} {{{\left| w \right|}^2}} dz \le
c{R^2}\int_{{B_{2R}}} {{{\left| {{D_0}w} \right|}^2}} dz +
c{R^2}\int_{{B_{2R}}} {\left( {{{\left| g \right|}^2} + {{\left| f
\right|}^2}} \right)} dz,
\end{equation}
it implies
\[\begin{array}{l}
\quad \int_{{B_R}} {{{\left| {{D_0}w} \right|}^2}} dz \\
  \le c\varepsilon {R^2}\int_{{B_{2R}}} {{{\left| {{D_0}w} \right|}^2}} dz + c\varepsilon {R^2}\int_{{B_{2R}}} {\left( {{{\left| g \right|}^2} + {{\left| f \right|}^2}} \right)} dz \\
 {\kern 1pt} {\kern 1pt} {\kern 1pt} {\kern 1pt} {\kern 1pt} {\kern 1pt} {\kern 1pt} {\kern 1pt}  + {c_\varepsilon }\int_{{B_R}} {\left( {{{\left| g \right|}^2} + {{\left| f \right|}^2}} \right)} dz + \varepsilon \int_{{B_R}} {{{\left| {{D_0}w} \right|}^2}} dz \\
  \le \varepsilon \int_{{B_{2R}}} {{{\left| {{D_0}w} \right|}^2}} dz + {c_\varepsilon }\int_{{B_{2R}}} {\left( {{{\left| g \right|}^2} + {{\left| f \right|}^2}} \right)} dz. \\
 \end{array}\]
Then for any $\rho  \le R$,
\[\begin{array}{l}
 \quad \int_{{B_\rho }} {{{\left| {{D_0}w} \right|}^2}} dz \le \int_{{B_R}} {{{\left| {{D_0}w} \right|}^2}} dz \\
  \le \varepsilon \int_{{B_{2R}}} {{{\left| {{D_0}w} \right|}^2}} dz + \frac{{{c_\varepsilon }{{\left( {2R - \rho } \right)}^2}}}{{{{\left( {2R - \rho } \right)}^2}}}\int_{{B_{2R}}} {{{\left| g \right|}^2}} dz + {c_\varepsilon }\int_{{B_{2R}}} {{{\left| f \right|}^2}} dz \\
  \le \varepsilon \int_{{B_{2R}}} {{{\left| {{D_0}w} \right|}^2}} dz + \frac{{{c_\varepsilon }{R^2}}}{{{{\left( {2R - \rho } \right)}^2}}}\int_{{B_{2R}}} {{{\left| g \right|}^2}} dz + {c_\varepsilon }\int_{{B_{2R}}} {{{\left| f \right|}^2}} dz. \\
 \end{array}\]
Now due to Lemma 2.9, we obtain
\[\int_{{B_\rho }} {{{\left| {{D_0}w} \right|}^2}} dz \le \frac{{c{R^2}}}{{{{\left( {2R - \rho } \right)}^2}}}\int_{{B_{2R}}} {{{\left| g \right|}^2}} dz + c\int_{{B_{2R}}} {{{\left| f \right|}^2}} dz,\]
and the conclusion holds with $\rho = R$.

\textbf{Theorem 6.2} Let $w \in W_{2,0}^{1,1}\left( \Omega \right)$
be a weak solution to (\ref{eq51}). Then for any $p \in \left[ {2,2
+ \frac{{2Q}}{{Q + 2}}{\varepsilon _0}} \right)$, we have ${D_0}w
\in L_{loc}^p(\Omega )$, and for any ${B_R} \subset {B_{4R}} \subset
\Omega $,
\begin{equation}
\label{eq56} \int_{{B_R}} {{{\left| {{D_0}w} \right|}^p}} dz \le
c\int_{{B_{4R}}} {\left( {{{\left| g \right|}^p} + {{\left| f
\right|}^p}} \right)} dz.
\end{equation}

\textbf{Proof}: By (4.2) and (\ref{eq52}), it follows
\[\begin{array}{l}
\quad \int_{{B_R}} {{{\left| {{D_0}w} \right|}^p}} dz \\
  \le c\left| {{B_R}} \right|{\left[ {{{\left( {\frac{1}{{\left| {{B_{2R}}} \right|}}\int_{{B_{2R}}} {{{\left| {{D_0}w} \right|}^2}} dz} \right)}^{\frac{1}{2}}} + {{\left( {\frac{1}{{\left| {{B_{2R}}} \right|}}\int_{{B_{2R}}} {{{\left( {{{\left| g \right|}^2} + {{\left| f \right|}^2}} \right)}^{\frac{p}{2}}}} dz} \right)}^{\frac{1}{p}}}} \right]^p} \\
  \le c\left| {{B_R}} \right|{\left[ {{{\left( {\frac{c}{{\left| {{B_{2R}}} \right|}}\int_{{B_{4R}}} {\left( {{{\left| g \right|}^2} + {{\left| f \right|}^2}} \right)} dz} \right)}^{\frac{1}{2}}} + {{\left( {\frac{1}{{\left| {{B_{2R}}} \right|}}\int_{{B_{2R}}} {{{\left( {{{\left| g \right|}^2} + {{\left| f \right|}^2}} \right)}^{\frac{p}{2}}}} dz} \right)}^{\frac{1}{p}}}} \right]^p} \\
  \le c\left| {{B_R}} \right|{\left[ {{{\left( {\frac{1}{{\left| {{B_{2R}}} \right|}}\int_{{B_{4R}}} {\left( {{{\left| g \right|}^p} + {{\left| f \right|}^p}} \right)} dz} \right)}^{\frac{1}{p}}} + {{\left( {\frac{1}{{\left| {{B_{2R}}} \right|}}\int_{{B_{2R}}} {\left( {{{\left| g \right|}^p} + {{\left| f \right|}^p}} \right)} dz} \right)}^{\frac{1}{p}}}} \right]^p} \\
  \le c\int_{{B_{4R}}} {\left( {{{\left| g \right|}^p} + {{\left| f \right|}^p}} \right)} dz. \\
 \end{array}\]

\textbf{Theorem 6.3 }Let $u \in W_2^{1,1}\left( \Omega  \right)$ be
a weak solution to (\ref{eq1}). Then for any $p \in \left[ {2,2 +
\frac{{2Q}}{{Q + 2}}{\varepsilon _0}} \right)$, we have ${D_0}u \in
L_{loc}^p(\Omega )$ and for any ${B_R} \subset {B_{4R}} \subset
\Omega$,
\begin{equation}
\label{eq57} \int_{{B_\rho }} {{{\left| {{D_0}u} \right|}^p}} dz \le
c\left[ {{{\left( {\frac{\rho }{R}} \right)}^{Q + 2 - \lambda
}}\int_{{B_{4R}}} {{{\left| {{D_0}u} \right|}^p}} dz + {\rho ^{Q + 2
- \lambda }}\left( {\left\| g \right\|_{{{L^{p,\lambda }}}}^p +
\left\| f \right\|_{{{L^{p,\lambda }}}}^p} \right)} \right].
\end{equation}

\textbf{Proof:} Combining Theorem 5.5 and Theorem 6.2 indicates
\begin{align}
&\quad \int_{{B_\rho }} {{{\left| {{D_0}u} \right|}^p}} dz
\le 2\int_{{B_\rho }} {{{\left| {{D_0}v} \right|}^p}} dz + 2\int_{{B_\rho }} {{{\left| {{D_0}w} \right|}^p}} dz \nonumber\\
  &\le c{\left( {\frac{\rho }{R}} \right)^{\frac{{2(Q + 2) - p(Q + 2 - \mu )}}{2}}}
  \int_{{B_R}} {{{\left| {{D_0}v} \right|}^p}} dz + 2\int_{{B_\rho }} {{{\left| {{D_0}w} \right|}^p}} dz \nonumber\\
 & \le c{\left( {\frac{\rho }{R}} \right)^{\frac{{2(Q + 2) - p(Q + 2 - \mu )}}{2}}}\int_{{B_R}} {{{\left| {{D_0}u} \right|}^p}} dz
 + c\int_{{B_R}} {{{\left| {{D_0}w} \right|}^p}} dz \nonumber\\
 & \le c{\left( {\frac{\rho }{R}} \right)^{\frac{{2(Q + 2) - p(Q + 2 - \mu )}}{2}}}
 \int_{{B_R}} {{{\left| {{D_0}u} \right|}^p}} dz + c\int_{{B_{4R}}}
 {\left( {{{\left| g \right|}^p} + {{\left| f \right|}^p}} \right)} dz\nonumber\\
  &\le c{\left( {\frac{\rho }{R}} \right)^{\frac{{2(Q + 2) - p(Q + 2 - \mu )}}{2}}}\int_{{B_R}} {{{\left| {{D_0}u} \right|}^p}} dz
  + c\frac{{\left| {{B_{4R}}} \right|}}{{{R^\lambda }}}\left( {\left\| g \right\|_{{{L^{p,\lambda }}}}^p
  + \left\| f \right\|_{{{L^{p,\lambda }}}}^p} \right) \nonumber\\
&\le c{\left( {\frac{\rho }{R}} \right)^{\frac{{2(Q + 2) - p(Q + 2 -
\mu )}}{2}}}\int_{{B_R}} {{{\left| {{D_0}u} \right|}^p}} dz + c{R^{Q
+ 2 - \lambda }}\left( {\left\| g \right\|_{{{L^{p,\lambda }}}}^p +
\left\| f \right\|_{{{L^{p,\lambda }}}}^p} \right).
\label{6.9}%
\end{align}
Let $H(\rho ) = \int_{{B_\rho }} {{{\left| {{D_0}u} \right|}^s}}
dz$, $H(R) = \int_{{B_R}} {{{\left| {{D_0}u} \right|}^s}} dz$,
${a_1} = \frac{{2(Q + 2) - s(Q + 2 - \mu )}}{2}$, ${b_1} = Q + 2 -
\lambda $, ${B_1} = c\left( {\left\| g \right\|_{L^{p,\lambda }}^p
+ \left\| f \right\|_{L^{p,\lambda }}^p} \right)$, $0 < \lambda  <
Q + 2$. Note that there exists $\mu$, $Q + 2 - \frac{{2\lambda }}{p} <
\mu < Q$ such that ${a_1} > {b_1}$. Hence we can conclude from Lemma 2.10
that
$$\int_{{B_\rho }} {{{\left| {{D_0}u} \right|}^p}} dz \le
c\left[ {{{\left( {\frac{\rho }{R}} \right)}^{Q + 2 - \lambda
}}\int_{{B_R}} {{{\left| {{D_0}u} \right|}^p}} dz + {\rho ^{Q + 2 -
\lambda }}\left( {\left\| g \right\|_{L^{p,\lambda }}^p + \left\|
f \right\|_{L^{p,\lambda }}^p} \right)} \right].$$

\textbf{Proof Theorem 1.2: } The result of Theorem 1.2 follows in virtue of Theorem 6.3 and the cutoff function technique.


{\bf Acknowledgements } This work is supported by the National
Natural Science Foundation of China(Grant Nos. 11271299, 11001221);
Natural Science Foundation Research Project of Shaanxi Province
(Grant No. JC201124).


\end{document}